%% file: Beam_Buckling_Manuscript/main.tex
\documentclass[11pt,a4paper]{article}

\title{Buckling Prediction for Nonlinear Elastic Beams with Soft Inclusions }

\author{%
Thien Tran-Duc%
\thanks{School of Mathematical Sciences, Adelaide University,  South Australia 5005.}
\footnote{\protect\url{mailto:thien.tran@adelaide.edu.au},
\protect\url{orcid:0000-0002-2004-5156}}
\and 
J.E. Bunder\footnotemark[1]
\footnote{\protect\url{mailto:Judith.Bunder@adelaide.edu.au},
\protect\url{orcid:0000-0001-5355-2288}}
\and 
A.J. Roberts\footnotemark[1]
\footnote{\protect\url{mailto:ProfAJRoberts@protonmail.com},
\protect\url{orcid:0000-0001-8930-1552}}
}

\date{\today}

\input{preamble}

\def\patRat{sparsity}% singular
\def\patRats{sparsities}% plural

\begin{document}

\maketitle

\begin{abstract}
%\todo{\par* What was done?
%\par* Why do it?
%\par* What were the results?
%\par* What do the results mean in theory and/or practise?
%\par* What is the reader's benefit?
%\par* How can readers use this information for themselves? 
%}
We develop an efficient and accurate multiscale computational framework for predicting the buckling and post-buckling behaviour of elastic beams containing periodically distributed soft inclusions. 
The framework extends our previous multiscale, patch, computational homogenisation for linear elasticity by incorporating nonlinearity, and thus enables accurate prediction of both the buckling onset and the subsequent post-buckling response. 
Microscale computations are performed only within a sparse set of small subdomains (patches), while the macroscale behaviour is recovered through a proven patch-coupling algorithm.
The scheme is assessed through quarter- and half-domain patch computations for beams with inclusion-to-matrix Young's modulus ratios ranging from~\(0.001\) to~\(1\). 
The results show that reducing the inclusion stiffness lowers both the critical buckling strain and the critical buckling stress, indicating an increased susceptibility to instability, while producing a milder post-buckling response with smaller transverse deflections and stress drops. 
Eigenvalue analysis of the Jacobian matrix accurately predicts the onset of instability and the corresponding critical strain and stress. 
Bifurcation diagrams of the nonlinear buckled configurations under compressive loading, and a quantitative analysis of the effect of the interpolation order on the predicted buckling and post-buckling responses, are also presented.
Comparisons with full-domain simulations demonstrate that the proposed framework accurately predicts both the buckling threshold and the post-buckling behaviour while substantially reducing the computational cost. 
The methodology is readily extendable to heterogeneous beams, plates, shells, and \text{other engineering structures.}
\end{abstract}

\tableofcontents

\section{Introduction}
\label{Sec: Intro}

\begin{figure}
\centering
\caption{\label{fig:Snap_1_15}Snapshots of snap-through buckling of a heterogeneous beam accurately predicted by computation on only the shown patches. 
The beam buckles under a uniform transverse load~\(F^e\) applied along the top edge of the central patch. 
Colours denote the axial normal stress~\(\sigma^{xx}\).
Only a quarter of the beam domain is computed~(\patRat~\cref{Epatrat} \(r=0.25\)).
}
\begin{enumerate}[leftmargin=0pt,nosep,label=(\alph*)]
\item[]\includegraphics[scale=0.8]{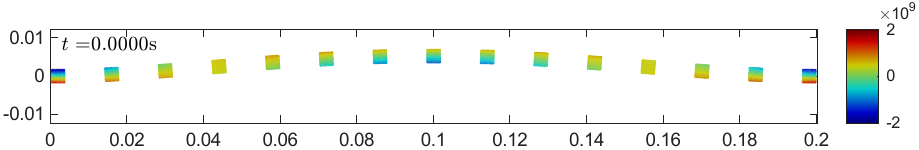}   
\item[]\includegraphics[scale=0.8]{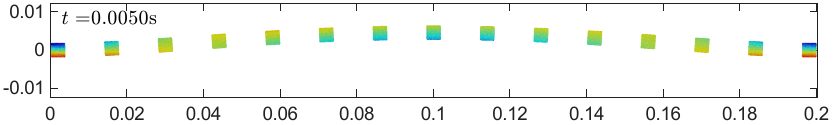}   
\item[]\includegraphics[scale=0.8]{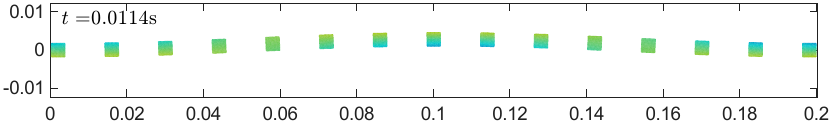}    
\item[]\includegraphics[scale=0.8]{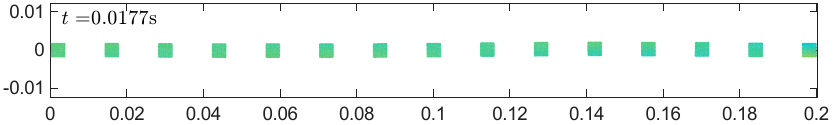}    
\item[]\includegraphics[scale=0.8]{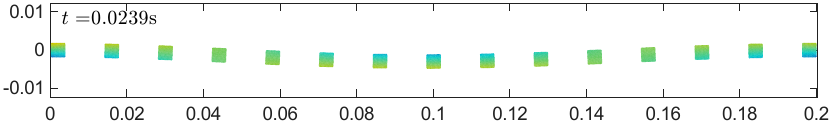}   
\item[]\includegraphics[scale=0.8]{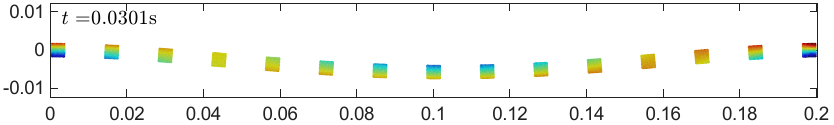}
\item[]\hfil space~\(x\)~(m)
\end{enumerate}
\end{figure}

Buckling instability in slender structures subjected to compressive loading is a phenomenon in which the structure undergoes a sudden and significant deformation, leading to a substantial reduction in load-carrying capacity and, in some cases, catastrophic failure~\citep{Han2013, Babilio2023, Riccio2025, Wagner2025}. 
Accurate prediction of buckling in mechanical structures poses a multiscale computational challenge because microscale material and geometric features may influence the onset and subsequent development of macroscale structural instability. 
Equation-free multiscale methods provide an efficient approach to such problems by performing detailed microscale simulations only within selected regions of the computational domain and appropriately coupling these simulations to recover system-level behaviour~\citep{Kevrekidis2003, Kevrekidis09a}. 
In particular, gap-tooth and related patch schemes exploit sparse microscale simulations to approximate macroscale dynamics without requiring the entire spatial domain to be microscopically resolved~\citep{Samaey03b, Hyman2005, Liu2015, Bunder2020a, Maclean2021}.
Related micro-to-macro computational strategies include \textsc{fe}\(^2\) methods, which couple microscale boundary-value problems to macroscale finite-element calculations~\citep{Raju2021, VBCTan2020}. 
The sparse simulation of nonlinear snap-through buckling shown in~\cref{fig:Snap_1_15} illustrates the capability of the patch approach that motivates its application here to heterogeneous beam buckling.

Buckling is a major concern in the design and operation of many engineering systems, including columns and beams in civil infrastructure~\citep{Karlsson2000}, thin-walled components in aerospace structures~\citep{Zhang2023}, offshore pipelines subjected to thermal and mechanical loads~\citep{Seth2021, Kyria2023}, ship hulls and marine structures~\citep{Yao2016, Liu2022}, and micro- and nano-scale devices~\citep{Xia2010, Hu2015}. 
Buckling occurs when the applied compressive load exceeds a critical value, beyond which the initial equilibrium configuration loses stability.
At this point, even small disturbances, geometric imperfections, or microscale heterogeneities may trigger a transition to a transversely deflected configuration, driven by the release of elastic energy stored under compression. 
The resulting instability can substantially reduce the load-carrying capacity and compromise the safety, reliability, and performance of engineering structures~\citep{Patra2023, Jia2024, Gulfam2025}.
Buckling of mechanical structures is inherently multiscale and highly nonlinear~\citep{Nezamabadi2010, Saavedra2012}, as microscale imperfections or material heterogeneities may influence a global structural instability. 
While a structure remains sufficiently far from its critical state, small geometric imperfections, local material heterogeneities, residual stresses, or slight variations in loading and boundary conditions may have only a minor influence on its mechanical response. 
However, as the applied load approaches a critical threshold, these perturbations may be rapidly amplified, destabilising the current equilibrium configuration and driving the structure towards a new equilibrium state. 
This sensitivity to small perturbations is a fundamental characteristic of buckling and presents a particular challenge for heterogeneous structures, where local microscale features may initiate or modify the macroscale instability. 

Homogenisation is a widely used multiscale technique for determining the effective mechanical properties of microstructured heterogeneous materials. 
It replaces the underlying complicated heterogeneous microstructure with an equivalent homogeneous continuum~\citep{Anthoine2010, Forest2011, Roberts2025a}. 
By averaging the microscale response, homogenisation significantly reduces computational cost while retaining the dominant macroscale behaviour of the material. 
Consequently, homogenised constitutive models have been extensively employed in analyses of buckling and stability of heterogeneous structures, where effective elastic properties are used to predict critical buckling loads, deformation modes, and post-buckling responses~\citep{ReiGon2016, Falsonea2021, Qi2025}. 
Such approaches have proved successful for many engineering materials with periodic or statistically homogeneous microstructures, particularly when the characteristic length scale of the heterogeneity is much smaller than that of the structure.
However, because microscale information is conventionally averaged, conventional homogenisation may not accurately capture localised deformation, stress concentrations, or instability mechanisms arising from strong material heterogeneity, defects, or non-periodic microstructures~\citep{Wagner2025b}. 
These limitations become increasingly important when local microscale features interact with \text{global structural instability.}

Computational multiscale methods seek to retain the relevant microscale information while reducing the cost associated with resolving the entire heterogeneous structure. 
Conventional full-domain simulations require the microstructure to be resolved throughout the computational domain and may therefore involve prohibitively large numbers of degrees of freedom for realistic engineering structures. 
\cref{fig:Snap_1_15} indicates how the Equation-Free gap-tooth and patch methodologies address this challenge by predicting macroscale behaviour from spatially sparse microscale simulations~\citep{Samaey03b, Hyman2005, Roberts2011a, Liu2015, Maclean2021}. 
By avoiding detailed microscale simulation over substantial portions of the physical domain, these approaches provide an efficient strategy for micro-to-macro computation. 
Related multiscale approaches include \textsc{fe}\(^2\) methods, which employ nested macro- and microscale finite-element computations \citep{Raju2021, VBCTan2020} and have also been applied to multiscale structural instability and buckling~\citep{Nezamabadi2010}.

The multiscale patch scheme employed herein \citep{Roberts2011a, Maclean2021, Bunder2021b} is a controllably accurate framework for predicting macroscale behaviour directly from microscale physics while substantially reducing the spatial domain that must be simulated (\cref{fig:Snap_1_15}). 
Originating from equation-free multiscale methodology and underpinned by  nonlinear dynamical systems theory \citep{Roberts1997, Kevrekidis2003, Roberts2007, Bunder2020a}, the scheme performs detailed microscale simulations only within a small number of disjoint patches distributed sparsely across the computational domain. 
The unsimulated space between neighbouring patches is bridged through interpolation and proven coupling conditions, see~\cref{subSec:PatchCoup}, enabling the macroscale solution to be accurately recovered \citep{Bunder2013b, Bunder2020a}. 
Importantly for buckling problems, the microscale model is retained within each patch, so that local material heterogeneities and their associated deformation are explicitly resolved rather than replaced by assumed effective material properties. 
Coupling between patches then communicates sufficient microscale information across the computational domain for the locally resolved physics to inform the global structural response.
Previous research developed this patch scheme for simulations of heterogeneous and functionally graded elastic beams \citep{TranD2024, TranD2025}. 
The scheme accurately predicted static equilibrium configurations and transient dynamic responses while simulating only a small fraction of the full computational domain. 
It captured the influence of spatially varying material properties on the deformation and stress fields of heterogeneous beams, with excellent agreement between patch and full-domain computations over a wide range of loading conditions. 
These studies established the patch scheme as an efficient and accurate multiscale computational framework for structural analysis. 
However, its capability to capture the onset of structural instability and the subsequent nonlinear post-buckling evolution of heterogeneous structures has \text{not yet been established.}

Herein, the patch methodology is extended to investigate the buckling and post-buckling behaviour of \emph{nonlinear} heterogeneous beams containing soft inclusions uniformly distributed in a single layer along the beam axis~(\cref{fig:BeamCoord}). 
The present study validates that the sparse patch simulations accurately capture the loss of stability and subsequent nonlinear evolution while retaining the microscale material heterogeneity within the simulated patches. 
\cref{Sec:NumRes} demonstrates that the patch methodology accurately predicts the buckling onset and post-buckling behaviour, including bifurcation points, equilibrium branches, and their associated stability characteristics, while retaining the computational advantage of resolving only a fraction of the \text{full structural domain.}

Although the present study specifically addresses beams with soft inclusions, the patch methodology is applicable to a broader range of micro-heterogeneous beams. 
For example, \cref{fig:Snap_1_15} illustrates a snap-through buckling process for a heterogeneous beam, with colour indicating the axial normal stress. 
The beam has length \(L=20\)\,cm and thickness \(W=0.35\)\,cm. 
Material heterogeneity is introduced through a spatially random and isotropic Young's modulus throughout the beam domain: \(E(x,y)\)~takes on iid random values distributed uniformly over \(180\,\text{GPa}<E(x,y)<200\)\,GPa. 
A uniformly distributed transverse load of \(F^e:=28.5\)\,MN/m is applied downward along the top edge of the central patch. 
Under this loading, the beam rapidly transitions from the initially upward-buckled configuration to a downward-buckled configuration via a snap-through event, with the transition occurring within approximately \(0.03\)\,s. 
Following snap-through, the beam remains buckled downwards while exhibiting decaying small-amplitude \text{oscillations due to inertia.}

\section{Patch algorithm for beam buckling}
\label{Sec:ProbDef}

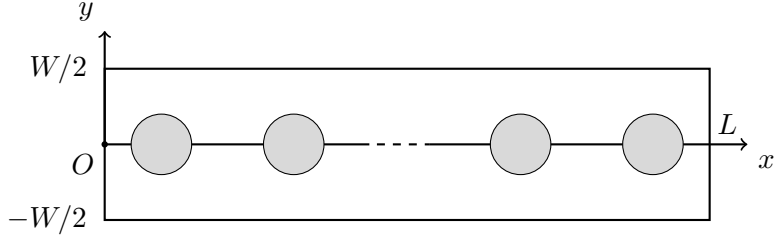
\begin{figure}
	\centering
	\caption{\label{fig:BeamCoord}
	Coordinate system for a two-dimensional beam with axial length~\(L\) and transverse width~\(W\), containing a single layer of identical circular soft inclusions distributed along the beam axis~(\(y=0\)).  }
	\input{Figs/BeamCoord.tex}    
\end{figure}

\cref{fig:BeamCoord} illustrates a 2-D beam with soft inclusions. 
The beam has length~\(L\) and width~\(W\). 
We place the origin at the centre of the left end of the beam, with the \(x\)-axis directed along the beam axis toward the right end, and the~\(y\)-axis oriented from the bottom edge to the top edge. 
 Accordingly, the undeformed beam's domain is~\((x,y)\in [0,L]\times[-W/2,W/2]\). 
All circular inclusions have the same radius~\(R_s\) and for a given beam, have the same Young's modulus that is less than that of the surrounding matrix material.
%Young's modulus is a measure of stiffness, defined as the ratio of compressive axial stress to axial strain, that is, the ratio of the compressive force per unit area to the spatial distortion per unit length. 
The inclusions with lower Young's moduli are characterised as `soft', or less stiff, as they are more easily distorted in response to an applied force.
Conversely, when  experiencing a specific distortion, the resulting stresses in the inclusions are smaller.
\cref{subSec:SubPatchDis,SSdiss} describe the material, microscale, spatial discretisation, while~\cref{SSpatches,subSec:PatchCoup} detail the scale-bridging patch algorithm for the accurate and efficient numerical exploration of buckling and related macroscale phenomena \text{in such elastic beams.}

\paragraph{Nonlinear elastic stress}
The 2-D beam is compressed along its axis, and thence many buckle. 
For finite deformations, the Green--Lagrange strain tensor is defined by the nonlinear relation that~\citep{Novoz1953}
%\todo{cite source of this nonlinear equation}
\begin{equation}
    {\varepsilonv}
    := \frac{1}{2}\left[
    \left(\nabla \uv \right)^T
    + \nabla \uv 
    + \left(\nabla \uv \right)^T\!\cdot\nabla \uv 
    \right],
    \label{eq:LagrGreenST}
\end{equation}
where~\(\uv :=(u,v)\) denotes the displacement vector, and the deformation gradient is evaluated with respect to the initial undeformed configuration. 
In two dimensions, the strain components are \text{explicitly written as~\citep{Orthwein1968}}
%\todo{cite source.}
\begin{equation}
\begin{aligned}
    \varepsilon^{xx}
    &= \frac{\partial u}{\partial x}
    + \frac{1}{2}\left[
    \left(\frac{\partial u}{\partial x}\right)^2
    + \left(\frac{\partial v}{\partial x}\right)^2
    \right], \\
    \varepsilon^{yy}
    &= \frac{\partial v}{\partial y}
    + \frac{1}{2}\left[
    \left(\frac{\partial u}{\partial y}\right)^2
    + \left(\frac{\partial v}{\partial y}\right)^2
    \right], \\
    \varepsilon^{xy}=\varepsilon^{yx}
    &= \frac{1}{2}\left(
    \frac{\partial u}{\partial y}
    + \frac{\partial v}{\partial x}
    + \frac{\partial u}{\partial x}\frac{\partial u}{\partial y}
    + \frac{\partial v}{\partial x}\frac{\partial v}{\partial y}
    \right).    
\end{aligned}\label{eq:STExplicitForm}
\end{equation}
Superscripts denote stress components~(\cref{subSec:SubPatchDis} introduces subscripts to denote grid indices in a microscale spatial discretisation).
The corresponding elastic stress induced by the finite deformation~\cref{eq:LagrGreenST} is described by the Saint--Venant--Kirchhoff model
\begin{equation}
    \sigmav 
    = \lambda \opn{Tr}(\varepsilonv){I}
    + 2\mu \varepsilonv,\label{eq:StressStrain}
\end{equation}
where \(\lambda(x,y)\) and~\(\mu(x,y)\) are the Lam\'e constants. 
%\todo{cite an example textbook source.}
This constitutive law corresponds to the strain-energy density 
\(\frac{\lambda}{2}\left[\opn{Tr}(\varepsilonv)\right]^2
    + \mu\,\opn{Tr}(\varepsilonv^2)\) \citep{Cihan2024}.
For other classes of hyperelastic materials, alternative strain-energy density functions may be adopted, leading to different constitutive relations~\cite[e.g.,][]{Orthwein1968, Guo2025}.
In component form, the stress tensor is
\begin{equation}
\begin{aligned}
    \sigma^{xx} &= (\lambda+2\mu)\varepsilon^{xx}  + \lambda\varepsilon^{yy},\\
    \sigma^{yy} &= \lambda\varepsilon^{xx}  + (\lambda+2\mu)\varepsilon^{yy},\\
    \sigma^{xy}&=\sigma^{yx} = 2\mu\varepsilon^{xy}.
\end{aligned}
\label{eq:StressComp}
\end{equation}

The elastic stress tensor is evaluated at material points throughout the beam. 
Spatial variation of the stress generates stress gradients and, consequently, internal forces that drive the motion of the material.
By Newton's second law, the acceleration of a material point satisfies \text{the \pde~\citep{Malvern1969}}
%\todo{cite a textbook source}
\begin{equation}
    \ddot{\uv }
    = \frac{1}{\rho}\left( \nabla \cdot\sigmav  + \Fv^e \right),
    \label{eq:NewSL}
\end{equation}
where~\(\Fv^e(x,y,t)\) denotes an applied external distributed body-force, if present, and~\(\rho(x,y)\) is \text{the material density.}

\paragraph{Boundary conditions}
To generate axial compression, asymmetric displacement boundary conditions are imposed at the two ends of the beam:
\begin{equation}
\uv (0,y) = (u^e,0), \qquad \uv (L,y) = (-u^e,0),
\label{eq:us}
\end{equation}
where \(u^e\) denotes the imposed horizontal compressive displacement. 
The corresponding axial strain is~\(\varepsilon^e={2u^e}/{L}\).

Along the top and bottom surfaces of the beam, stress-free boundary conditions are imposed:~\(\sigmav\cdot\nv  = 0\), or in component form
\begin{equation}
	\sigma^{xx}n^x + \sigma^{xy}n^y = 0,\qquad
	\sigma^{xy}n^x + \sigma^{yy}n^y = 0,
	\label{eq:SFreeExpl}
\end{equation}
where~\(\nv :=(n^x,n^y)\) denotes the outward unit normal vector at a point~\((x,y)\) \text{on the boundary.} 

\paragraph{Non-dimensionalization} Using non-dimensional variables 
\(\uv^*=\uv/L_0\), \(t^*=t/t_0\), \(\boldsymbol{\nabla}^*=L_0\boldsymbol{\nabla}\), 
\(\sigma^*=\sigma/\sigma_0\), and \((\Fv^e)^*=\Fv^e/\Fv^e_0\), 
and dividing the governing Newton's law \pde~\cref{eq:NewSL} by~\(\rho L_0/t_0^2\), its non-dimensional form is
\begin{equation*} 
\ddot{\uv}^*=\frac{\sigma_0t_0^2}{L_0^2\rho}\boldsymbol{\nabla}^* \cdot\sigma^* +\frac{\Fv^e_0t_0^2}{L_0\rho}(\Fv ^e)^* .%\label{eqn:NonDimDisp}
\end{equation*}
We choose the characteristic scales 
\(L_0 = W\), 
\(t_0 = L_0\sqrt{\rho_m/E_m}\), 
\(\sigma_0 = L_0^2\rho_m/t_0^2\), 
and \(\Fv^e_0 = L_0\rho_m/t_0^2\), 
where~\(W\), \(\rho_m\) and~\(E_m\) denote, respectively, the beam width, density and Young's modulus of the matrix material.
The non-dimensional equation then simplifies to
\begin{equation}
    \ddot{\uv}^* = \frac{1}{\rho^*}\left(\boldsymbol{\nabla}^* \cdot \sigma^* + (\Fv^e)^*\right),
    \label{Eqn:StressDiv}
\end{equation}
in which~\(\rho^*\le 1\) is the dimensionless (relative) density.
Hereafter, the superscript~\(*\) is omitted for notational simplicity, and then~\cref{Eqn:StressDiv} takes the same form \text{as that of~\cref{eq:NewSL}.}

\subsection{Patches}
\label{SSpatches}

\begin{figure}
\centering
\caption{\label{fig:PatchLayout}
Layout of the patches in the patch scheme for different patch coverage levels: (a)~\(I_p=I_L\)  patches (full domain, \(r=1\)), (b)~\(I_p=\left(I_L+1\right)/2\) patches (half domain, \(r=0.5\)), (c)~\(I_p=\left(I_L+3\right)/4\) patches (quarter domain, \(r=0.25\)). }
\input{Figs/PatchLayout}
\end{figure}
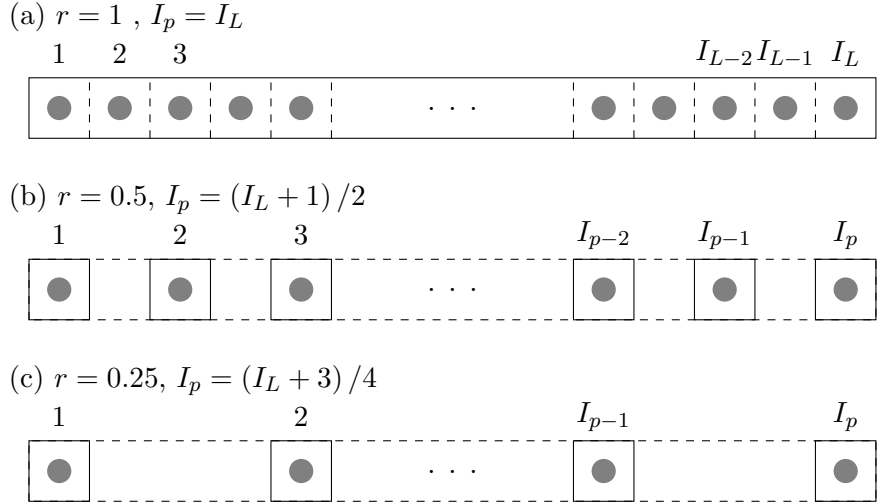%
The key components of the patch scheme~\cite[e.g.,][]{Kevrekidis09a, Bunder2020a} are the definition, distribution, and coupling of the patches.  
For a \emph{full-domain} computation, the beam domain is partitioned along its axis into~\(I_L\) equal size subdomains, called patches, as illustrated in~\cref{fig:PatchLayout}(a).  
Each patch contains a single soft inclusion located at its centre and has length~\(h:=L/I_L\).
The number of patches employed in a patch scheme prediction is some chosen number~\(I_p\ll I_L\).  
When~\(I_p=I_L\), the entire beam domain is fully covered by patches: adjacent patches share \(x\)-edges: such a full-domain computation is equivalent to a fully resolved microscale simulation over the entire beam domain.
For~\(I_p<I_L\), the patches are sparsely distributed and separated by a uniform centre-to-centre distance~\(H\).  
The ratio of the patch size~\(h\) to the patch spacing~\(H\), termed the \emph{\patRat}, \text{is defined to be}
\begin{equation}
r:=\frac{h}{H}\le 1\,.  \label{Epatrat}
\end{equation}
Smaller values of \patRat~\(r\) correspond to smaller domain coverage and hence greater computational efficiency.
It is then the patch spacing~\(H\) that dominantly controls the macroscale accuracy, as discussed in~\cref{Sec:NumRes}.
Examples of patch layouts corresponding to different coverage levels are shown in~\cref{fig:PatchLayout}(a)--(c), where~\(I_p=I_L\), \(I_p=0.5\left(I_L+1\right)\), and~\(I_p=0.25\left(I_L+3\right)\), yielding \patRats\ \(r=1\), \(0.5\), and~\(0.25\), respectively.
Within each patch, the microscale elastic model described above is solved independently~(\cref{SSdiss,subSec:SubPatchDis}) and their \(x\)-edge values are set via the inter-patch interpolation of~\cref{subSec:PatchCoup} that couples the patches \text{into a macroscale system.}

\subsection{Microscale spatial discretisation within patches} 
\label{subSec:SubPatchDis}

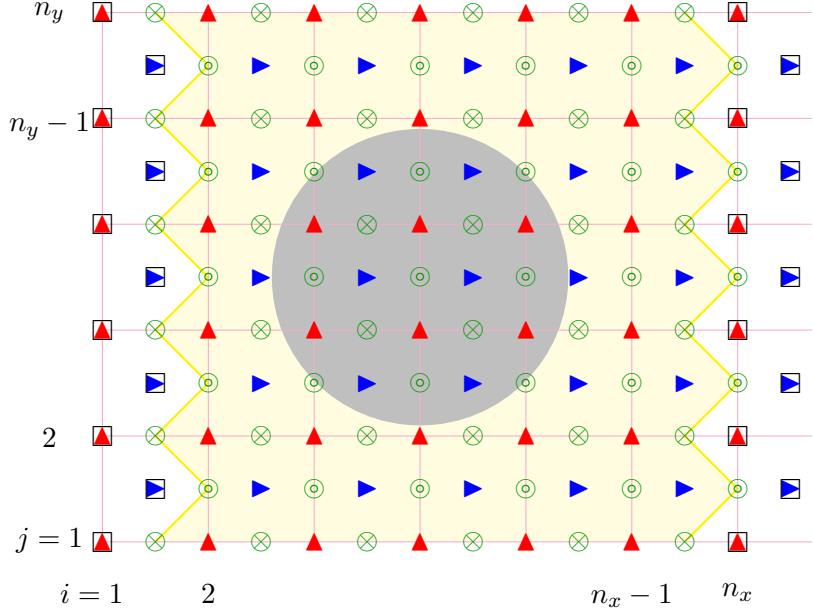
\begin{figure}
\centering
\caption{\label{fig:Microgrid}%
Microscale staggered \(xy\)-grid in a representative patch (shaded yellow matrix, and gray inclusion) of the 2D beam (magenta grid denotes integer indices~\(i,j\)). 
The grid symbols indicate: \uSym~horizontal displacement~\(u\); \vSym~vertical displacement~\(v\); 
\xSym~shear strains~\(\varepsilon^{xy}=\varepsilon^{yx}\) and shear stresses~\(\sigma^{xy}, \sigma^{yx}\); 
\oSym~normal strains~\(\varepsilon^{xx}, \varepsilon^{yy}\) and normal stresses~\(\sigma^{xx}, \sigma^{yy}\); $\square$ patch-interpolated values.
The gray solid circle represents the soft inclusion.}
\input{Figs/Grid.tex}
\end{figure}

The equations of motion~\cref{Eqn:StressDiv} are numerically solved within each patch using sub-patch, centred, finite difference, schemes on a staggered micro-grid. 
\cref{fig:Microgrid} illustrates that in this staggered arrangement the displacement and stress fields are evaluated at different micro-spatial locations to enhance numerical stability and accuracy. 
A two-dimensional Cartesian grid of size~\(n_x \times n_y\) is used to discretise each patch, with grid spacings~\(\delta x\) and~\(\delta y\) in the horizontal and vertical directions, respectively. 
The leftmost and rightmost vertical grid layers (unshaded in~\cref{fig:Microgrid}) are located on the edges of the patch domain, and serve as \text{the \(x\)-boundaries of each patch.}

The field variables are defined on the staggered grid as follows, where \emph{cell} refers to each sub-patch micro-grid square (thin magenta lines in \cref{fig:Microgrid}).
\begin{itemize}
    \item The horizontal displacement~\(u\) is evaluated at cell centres, that is at half-integer locations in both directions,~\(u_{i+1/2,\,j+1/2}\).
    \item The vertical displacement~\(v\) is evaluated at grid nodes, that is at integer locations,~\(v_{i,j}\).
    \item The normal stresses~\(\sigma^{xx}\) and~\(\sigma^{yy}\) are evaluated at the centres of cell vertical-edges,  \(\sigma^{xx}_{i,j+1/2}\) and~\(\sigma^{yy}_{i,j+1/2}\).
    \item The shear stress~\(\sigma^{xy}\) is evaluated at the centres of cell horizontal-edges,~\(\sigma^{xy}_{i+1/2,j}\)\,.
\end{itemize}
Accordingly, the strain components~\cref{eq:STExplicitForm} in any patch, \(I\in\{1,\ldots, I_p\}\), are discretised using centred finite differences on the staggered grid:
 \begin{equation}
\begin{aligned}
\varepsilon^{xx,I}_{i,j\phalf} &= \frac{\delta_{i}u^I_{i,j\phalf}}{\delta x}+\frac{1}{2}\left[\left(\frac{\delta_{i}u^I_{i,j\phalf}}{\delta x}\right)^2+\left(\frac{\delta_{i}v^I_{i,j\phalf}}{\delta x}\right)^2\right],\\
\varepsilon^{yy,I}_{i,j\phalf} &= \frac{\delta_{j}v^I_{i,j\phalf}}{\delta y}+\frac{1}{2}\left[\left(\frac{\delta_{j}u^I_{i,j\phalf}}{\delta y}\right)^2+\left(\frac{\delta_{j}v^I_{i,j\phalf}}{\delta y}\right)^2\right],\\
\varepsilon^{xy,I}_{i\phalf,j} &= \frac{1}{2}\left(\frac{\delta_{j}u^I_{i\phalf,j}}{\delta y}+\frac{\delta_{i}v^I_{i\phalf,j}}{\delta x}+\frac{\delta_{i}u^I_{i\phalf,j}}{\delta x}\frac{\delta_{j}u^I_{i\phalf,j}}{\delta y}+\frac{\delta_{i}v^I_{i\phalf,j}}{\delta x}\frac{\delta_{j}v^I_{i\phalf,j}}{\delta y}\right).
\end{aligned}\label{eq:Dis_Strain}
\end{equation}
The symbols~\(\delta_i\) and~\(\delta_j\) denote centred finite differences of a field~\(f_{i,j}\) in the \(x\)- and \(y\)-directions, respectively:  \(\delta_i f_{i,j} := f_{i\phalf,j} - f_{i\mhalf,j}\) and \(\delta_j f_{i,j} := f_{i,j\phalf} - f_{i,j\mhalf}\).
Here and hereafter, the index superscripts~\(\phalf\) and~\(\mhalf\) respectively denote shifts by positive and negative half-integer indices.
On the staggered grid~(\cref{fig:Microgrid}), the horizontal displacement \(u^I_{i\phalf,\,j\phalf}\) and the vertical displacement~\(v^I_{i,j}\) are defined at half-integer locations~\((i\phalf,j\phalf)\) and integer locations~\((i,j)\), respectively. 
However, the finite difference operators~\(\delta_j u^I_{i,j\phalf}\), \(\delta_i u^I_{i\phalf,j}\), \(\delta_i v^I_{i,j\phalf}\), and~\(\delta_j v^I_{i\phalf,j}\) require values of~\(u\) and~\(v\) at integer and half-integer locations, respectively. 
To compute the strain~\cref{eq:Dis_Strain}, obtain these values by interpolating from neighbouring grid points: 
\(u^I_{i,j} = \tfrac{1}{4}\big(u^I_{i\phalf,j\phalf}+u^I_{i\mhalf,j\phalf}+u^I_{i\phalf,j\mhalf}+u^I_{i\mhalf,j\mhalf}\big)\), and 
\(v^I_{i\phalf,j\phalf} = \tfrac{1}{4}\big(v^I_{i+1,j}+v^I_{i,j,N}+v^I_{i,j+1}+v^I_{i+1,j+1}\big)\).

Once the strain components~\cref{eq:Dis_Strain} are evaluated, the corresponding stress components~\cref{eq:StressComp} are computed using
\begin{equation}
	\begin{aligned}
		\sigma^{xx,I}_{i,j\phalf} &= (\lambda^I_{i,j\phalf}+2\mu^I_{i,j\phalf})\varepsilon^{xx,I}_{i,j\phalf}  + \lambda^I_{i,j\phalf}\varepsilon^{yy,I}_{i,j\phalf},\\
		\sigma^{yy,I}_{i,j\phalf} &= \lambda^I_{i,j\phalf}\varepsilon^{xx,I}_{i,j\phalf}  + (\lambda^I_{i,j\phalf}+2\mu^I_{i,j\phalf})\varepsilon^{yy,I}_{i,j\phalf},\\
		\sigma^{xy,I}_{i\phalf,j} &=\sigma^{yx,I}_{i\phalf,j,N}  = 2\mu^I_{i\phalf,j}\varepsilon^{xy,I}_{i\phalf,j}.
	\end{aligned}\label{eq:Dis_Stress}
\end{equation} 
Once the stress components~\cref{eq:Dis_Stress,eq:sxysyysxx} have been evaluated, the acceleration components in Newton's law \pde~\cref{Eqn:StressDiv} are computed. 
Owing to the staggered grid arrangement, the acceleration components are evaluated by centred differences on the micro-grid of~\cref{fig:Microgrid}:
\begin{subequations}\label{eq:DisGovE}%
\begin{align}
\ddot{u}^I_{i\phalf,j\phalf}&=\frac{1}{\rho^{I}_{i\phalf,j\phalf}}\left[\left(\frac{\delta_i\sigma^{xx,I}_{i\phalf,j\phalf}}{\delta x}+\frac{\delta_i\sigma^{xy,I}_{i\phalf,j\phalf}}{\delta y}\right) + F^{x,I}_{i\phalf,j\phalf}\right], 
\label{EdisGovu}\\
\ddot{v}^I_{i,j} &= \frac{1}{\rho^I_{i,j}}\left[\left(\frac{\delta_i\sigma^{xy,I}_{i,j}}{\delta x}+\frac{\delta_i\sigma^{yy,I}_{i,j}}{\delta y}\right) + F^{y,I}_{i,j}\right]. 
\label{EdisGovv}
\end{align}
\end{subequations}
To close the computation we need the following boundary and edge values around \text{the patches of~\cref{fig:Microgrid}.}

\paragraph{Beam top and bottom boundary conditions}
Along the top and bottom boundaries of the beam, at~\(y=\pm W/2\) corresponding to~\(j=1\) and~\(j=n_y\), the stress-free boundary conditions~\cref{eq:SFreeExpl} become
\begin{equation}
\begin{aligned}
\sigma^{yy,I}_{i,n_y} = \left(\frac{n^{x,I}_{i,n_y}}{n^y_{i,n_y}}\right)^2\sigma^{xx,I}_{i,n_y}, & \qquad
\sigma^{yy,I}_{i,1}    = \left(\frac{n^{x,I}_{i,1}}{n^y_{i,1}}\right)^2\sigma^{xx,I}_{i,1},\\
\sigma^{xy,I}_{i,n_y} = -\frac{n^{x,I}_{i,n_y}}{n^y_{i,n_y}}\sigma^{xx,I}_{i,n_y}, &\qquad
\sigma^{xy,I}_{i,1}    = -\frac{n^{x,I}_{i,1,N}}{n^{y,I}_{i,1}}\sigma^{xx,I}_{i,1}.
\end{aligned}\label{eq:sxysyysxx}
\end{equation} 
Hence, first~\(\sigma^{xx,I}_{i,1}\) and~\(\sigma^{xx,I}_{i,n_y}\) are evaluated, and then they determine the remaining stress components.
Since~\(\sigma^{xx,I}\) depends only on the derivative of the displacement in the~\(x\)-direction, its evaluation requires displacement values only along the top and bottom boundaries. 
While the discrete vertical displacement~\(v\) is defined directly at grid points on these boundaries, the discrete horizontal displacement~\(u\) is not defined, due to the staggered grid arrangement. 
Therefore, \(u\) is linearly extrapolated from nearby micro-grid values:
\(u^I_{i\phalf,n_y} = 1.5u^I_{i\phalf,n_y\mhalf}-0.5u^I_{i\phalf,n_y-1\mhalf}\), and~\(u^I_{i\phalf,1} = 1.5u^I_{i\phalf,1\phalf}-0.5u^I_{i\phalf,2\phalf}\).
These extrapolated values are then used to evaluate~\(\sigma^{xx,I}\) at the top and \text{bottom beam boundaries.}

The stress equations~\cref{eq:Dis_Stress} determine the normal stresses~\(\sigma^{xx,I}_{i,1}\) and~\(\sigma^{xx,I}_{i,n_y}\).
Then~\cref{eq:sxysyysxx} calculates the corresponding boundary shear and normal stresses, \(\sigma^{xy,I}_{i,1}\), \(\sigma^{xy,I}_{i,n_y}\), \(\sigma^{yy,I}_{i,1}\), and~\(\sigma^{yy,I}_{i,n_y}\). 
The unit normal vector at the corresponding grid nodes along the boundaries is computed from the boundary displacement field as
\begin{equation*}
\begin{aligned}
n^{x,I}_{i,1}     &=\delta_iv^I_{i,1}, \qquad
n^{y,I}_{i,1}       =-\delta_iu^I_{i,1}-\delta x,\\
n^{x,I}_{i,n_y} &=-\delta_iv^I_{i,n_y}, \qquad 
n^{y,I}_{i,n_y}    =\delta_iu^I_{i,n_y}+\delta x.\\
\end{aligned}
\end{equation*}

\paragraph{Patch left-right edge conditions}
\cref{subSec:PatchCoup} describes the patch coupling interpolation that determines the left-edge displacement values, \(u^I_{i+1/2,\,j+1/2}\) and~\(v^I_{i,j}\) at~\(i=1\) for~\(I=2,\ldots,I_p\).
For the leftmost end patch~\(I=1\) we use the strain-controlled boundary condition~\cref{eq:us}:
\(u^1_{3/2,\,j+1/2}= u^{e}\); 
\(v^1_{1,j}= 0\).
Similarly, the patch coupling interpolation determines the right-edge displacement values, \(u^I_{n_x+1/2,\,j+1/2}\) and~\(v^I_{n_x,j}\) for~\(I=1,\ldots,I_p-1\). 
For the rightmost end patch~\(I=I_p\) we use the strain-controlled boundary condition~\cref{eq:us}: 
\(u^{I_p}_{n_x+1/2,\,j+1/2}= -u^{e}\); 
\(v^{I_p}_{n_x,j}= 0\).

\subsection{Patch coupling}
\label{subSec:PatchCoup}

The computation is performed on a set of patches that are sparsely distributed along the axial direction of the beam.%
\footnote{In beams with greater macro- to macro-scale separation one could achieve greater efficiencies by invoking sparse patches across the beam as well as along the beam.}
These patches are coupled across the unsimulated space to accurately predict the macroscale dynamics.
The adopted patch-coupling, developed by~\cite{Bunder2020a, Roberts2023a} and implemented in the Toolbox~\cite[]{Maclean2021}, preserves the self-adjoint symmetry of the underlying microscale model, thereby enhancing numerical stability in the macroscale simulations.
Amazingly, coupling via standard polynomial interpolation has been proved \text{to be good.}

We choose an order~\(m\) (even) for the inter-patch polynomial interpolation, denoted~\(P_m\).
Herein we explore~\(m\in\{2,4,6,8\}\).
Then left- and right-edge values for each patch~(\cref{fig:Microgrid}) are determined by interpolation from the right- and left-next-to-edge values of itself and its \(m\)~nearest neighbour patches.
That is, for the displacement fields~\((u,v)\), left edge values~(\(i=1\)) of each patch~\(I\) are interpolated using the right next-to-edge values~(\(i=n_x-1\)) from the patch~\(I\), together with corresponding values from \(m\)~neighbouring patches to the left and right:  away from the beam ends, the \(m\)~neighbours are indexed by~\(I\mhalf := \{I-m/2, \ldots, I-1\}\) and~\(I\phalf := \{I+1, \ldots, I+m/2\}\). 
Therefore:
\begin{subequations}\label{EEedgeInterp}%
\begin{equation}
\begin{aligned}
u^I_{1\phalf,j\phalf} &= P_m(u^I_{n_x-1\phalf,j\phalf}, u^{I\mhalf}_{n_x-1\phalf,j\phalf},u^{I\phalf}_{n_x-1\phalf,j\phalf}),\\
v^I_{1,j} &= P_m(u^I_{n_x-1,j}, u^{I\mhalf}_{n_x-1,j},u^{I\phalf}_{n_x-1,j}).
\end{aligned}\label{eq:LeftInterp}
\end{equation}
Similarly, the right edge values (\(i=n_x\)) of patch~\(I\) are interpolated using the left next-to-edge values (\(i=2\)) of patch~\(I\) and its~\(m/2\) left neighbouring patches and~\(m/2\) right neighbouring patches:
\begin{equation}
\begin{aligned}
u^I_{n_x\phalf,j\phalf} &= P_m(u^I_{2\phalf,j\phalf}, u^{I\mhalf}_{2\phalf,j\phalf},u^{I\phalf}_{2\phalf,j\phalf}),\\
v^I_{n_x,j} &= P_m(u^{I}_{2,j}, u^{I\mhalf}_{2,j},u^{I\phalf}_{2,j}).\\
\end{aligned}\label{eq:RighttInterp}
\end{equation}
\end{subequations}
These interpolated patch-edge values determine the macroscale feedback into the microscale sub-patch physical computations.

This patch scheme also encompasses full-domain computations when we craft the patches so that the right-edge of patch~\(I\) is at precisely the same~\(x\) as the left-next-to-edge of patch~\(I+1\).
Then the inter-patch interpolation, for every~\(m\), reduces to the trivial copy operation, and hence the patch scheme reduces to \text{a full-domain computation.}

The spacing between interpolated displacement fields is the patch spacing~\(H\).
This spacing is usually chosen by the desired accuracy at the macroscale of the scenario to be considered.
Earlier work~\cite[e.g.,][]{Bunder2020a} has proven and verified that in general the accuracy of the patch scheme with \(m\)th~order interpolation  varies with patch spacing~\(H\) \text{to an error~\Ord{H^m}.}

\paragraph{Patch size}
A key advantage of the patch scheme is that microscale computations are performed only within patches that occupy a small fraction of the overall domain. 
While smaller patches lead to greater computational savings, an important question is how small the patches can be chosen without compromising accuracy.
For materials with periodic or repeating microstructures, each patch should span at least one period of the microstructure. 
In such cases, a patch is analogous to a representative volume element used in homogenisation techniques. 
Although the patch size may also be a multiple of the microstructural period, this must be balanced against computational efficiency.
For highly heterogeneous materials, the patch size must be sufficiently large to capture the essential macroscale dynamics and to mitigate the influence of localised microscale fluctuations, while remaining small enough to retain computational efficiency.
Accordingly, as shown in~\cref{fig:Microgrid}, here each patch is chosen to contain one microscale period with its \text{soft inclusion centred}.

\subsection{Dissipation}
\label{SSdiss}

Weak dissipation may arise due to boundary friction with the ambient environment, acoustic radiation to the surroundings, and\slash or internal viscoelastic effects. 
To account qualitatively for such dissipative mechanisms we introduce a Kelvin--Voigt viscoelastic term, \(+\eta \nabla^2 \dot{\uv}\), into the non-dimensional \pde~\cref{Eqn:StressDiv}. 
This contribution is incorporated into the microscale discretisation~\cref{eq:DisGovE}: 
\begin{subequations}\label{EEvis}%
\begin{align}
\ddot{u}^I_{i\phalf,j\phalf}
&= \cref{EdisGovu} +
\eta\left(
\frac{\delta_i^2\dot{u}^I_{i\phalf,j\phalf}}{\delta x ^2}
+\frac{\delta_j^2\dot{u}^I_{i\phalf,j\phalf}}{\delta y ^2}
\right);
\label{eq:Visu}
\\
\ddot v^I_{i,j}
&= \cref{EdisGovv}+
\eta\left(
\frac{\delta_i^2\dot{v}^I_{i,j}}{\delta x ^2}
+\frac{\delta_j^2\dot{v}^I_{i,j}}{\delta y ^2}
\right),
\label{eq:Visv}
\end{align}
\end{subequations}
 in which the second-order finite differences are defined as $\delta_i^2f_{i,j}:=\delta_i\left(\delta_if_{i,j}\right)=f_{i+1,j}-2f_{i,j}+f_{i-1,j}$, and $\delta_j^2f_{i,j}:=\delta_j\left(\delta_jf_{i,j}\right)=f_{i,j+1}-2f_{i,j}+f_{i,j-1}$.
This Kelvin--Voigt viscoelasticity is also representative of a range of other phenomenological dissipation with non-dimensional strength~\(\eta\): herein we specify \text{weak dissipation with \(\eta := 10^{-3}\).}

\section{Multiscale predictions of buckling}
\label{Sec:NumRes}

For the purposes of validation, we compute and analyse both a homogeneous beam and various soft-inclusion beams of contrasting Young's moduli.  
Comparison reveals how the presence of soft inclusions alters the buckling characteristics \text{of a beam.}

All quantities in this section are non-dimensional.
The beam has length~\(L:=57\) and width~\(W:=1\).  
The Young's modulus and density of the matrix material are both set to~\(1\), and the Poisson ratio of the matrix material is \(\nu:=0.3\).  
The soft inclusions have radius \(R_s:=0.3\), Young's modulus~\(E_s:=0.1\), \(0.01\), or~\(0.001\), density~\(\rho_s:=1\), and Poisson ratio~\(\nu_s:=0.3\).  
The density of the soft inclusions does not affect the equilibrium states of the buckled beam.  
However, the density affects the beam dynamics, since the acceleration of material points is inversely proportional to \text{the local material density.}

The beam domain is partitioned into~\(57\) equal subdomains, or patches, each of length~\(h:=1\), with one soft inclusion located at the centre of each subdomain.  
Consequently, full-domain computation is performed over the entire beam domain using~\(I_p=57\) patches.  
This configuration serves as the benchmark for comparison with sparse patch computations using~\(I_p=15\) and~\(I_p=29\) patches, with  respective \patRats~\cref{Epatrat} \(r=0.25\) and~\(r=0.5\). 
Unless otherwise stated, our numerical results are obtained using a  sub-patch micro-grid with~\(n_x:=15\) and~\(n_y:=14\).
Here, \(n_y=n_x-1\) since the left and right interpolated grid layers of each patch are located on the left\slash right edges of each patch, as illustrated by the micro-grid in \cref{fig:Microgrid}.
The micro-grid determines the microscale spatial resolution, \text{here~\(\delta x=\delta y=0.077\).}

The numerical computations require the displacement vector~\(\uv\) and its time derivative~\(\dot{\uv}\), with horizontal and vertical components evaluated at sub-patch micro-grid locations, respectively:
\begin{align}&
\uv =
\left(u^I_{i+1/2,j+1/2},v^I_{i,j}\right),
&&
\dot{\uv} =
\left(\dot{u}^I_{i+1/2,j+1/2},\dot{v}^I_{i,j}\right)
\end{align}
for micro-grid~\(i=1,\ldots,n_x\), \(j=1,\ldots,n_y\), and patches~\(I=1,\ldots,I_p\).
At any time, the state of the beam is then described by the state vector~\(\Uv:=(\uv,\dot\uv)\).
The discretisation of the dynamic \pde~\cref{Eqn:StressDiv}, plus the dissipative terms~\cref{EEvis},  provide the dynamical equations for the beam: namely~\(\dot{\Uv} = \Fv(\Uv)\).
The nonlinear function~\(\Fv\) is defined by our coded discretisations of the stress field divergence~\(\nabla\cdot\sigmav\), the dissipation~\(\eta \nabla^2 \dot{\uv}\), and \text{any external force~\(\Fv^e\). }

The steady-state form of the evolution equation ~\(\dot{\Uv} = \Fv(\Uv)\), namely~\(\Fv(\Uv)=\vec0\), determines equilibrium configurations of the beams~(\cref{SSdsf}), denoted by~\(\Uv_e\).
The onset of buckling from an equilibrium compression~\(\Uv_e\) (\cref{subSec:Emodes}) is analysed through the eigenvalues and eigenvectors of the Jacobian matrix
\begin{equation}
\mathbf{J}:=\left.\D\Uv\Fv\right|_{\Uv=\Uv_e}.
\label{eq:JabMatrix}
\end{equation}

%\todo{Throughout this section cite any comparable work and results}

\subsection{Displacement and stress of compressed equilibria}
\label{SSdsf}

\begin{figure}
\centering
\caption{\label{fig:HomoStress}Axial normal stress~\(\sigma^{xx}\) (top row), and shear stress~\(\sigma^{xy}\) (bottom row), in the three leftmost patches of the homogeneous beam~(\(E_s=0.01\)) under axial compression strain~\(\varepsilon^e=0.4\%\), obtained from the full-domain computation~(\(r=1\)).
This strain compresses the beam to~\(0.11<x<56.89\)\,.
}
\begin{tabular}{@{}c@{}c@{}c@{}}
\includegraphics{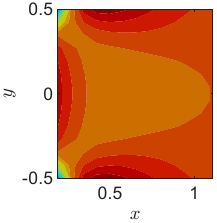} & 
\includegraphics{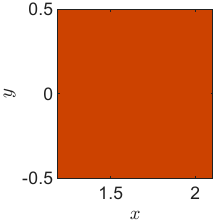} & 
\includegraphics{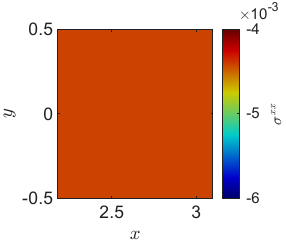}\\
\includegraphics{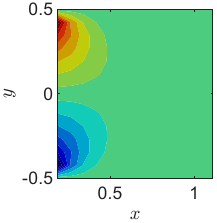} & 
\includegraphics{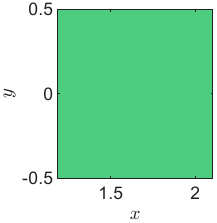} & 
\includegraphics{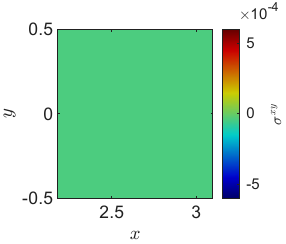}\\ 
\end{tabular}
\end{figure}%

Axial compression is applied by specifying horizontal displacements at the two ends of the beam via the boundary conditions~\cref{eq:us}.
With equal-and-opposite axial compression applied to both ends, the left\slash right half of the beam moves right\slash left towards the centre of the beam at~\(x=L/2\). 
Consequently, the centre of the beam~(\(x=28.5\)) is \text{stationary, unless buckled.}

For compressed axial strain~\(\varepsilon^e=0.4\%\)  in a homogeneous beam, \Cref{fig:HomoStress} shows the axial normal stress~\(\sigma^{xx}\) and the shear stress~\(\sigma^{xy}\) in the three leftmost patches. 
The applied compression generates non-uniform stress fields within narrow boundary layers which only extend a short distance into the beam~\citep{Timo1951}: here they extend less than the length of the leftmost patch.
Outside these boundary layers, \(\sigma^{xx}\) is constant and~\(\sigma^{xy}=0\),  as seen in the second and third leftmost patches.
The uniformity of~\(\sigma^{xx}\) outside the boundary layer corresponds to uniform normal strains~\(\varepsilon^{xx}\) and~\(\varepsilon^{yy}\), as implied by the linear constitutive relation~\cref{eq:StressStrain} and the constant elastic properties of the homogeneous material.
From the nonlinear strain equations~\cref{eq:STExplicitForm}, uniform normal strains imply a linear profile of both the horizontal displacement~\(u\) along the beam axis and  the vertical displacement~\(v\) across \text{the beam cross-section.}

\begin{figure}
\centering
\caption{\label{fig:uvtbedgesEs1001}
Displacement profiles and shapes of the top and bottom beam edges for the homogeneous beam~(\(E_s=1\), blue curves and symbols) and the soft-inclusion beam~(\(E_s=0.01\), red curves and symbols) under axial compression~\(\varepsilon^e=0.4\%\).
The dotted curves represent the full-domain computations~(\(r=1\)), and the~\(+\) and~\(\times\) symbols denote half-domain~(\(r=0.5\)) and quarter-domain~(\(r=0.25\)) computations, respectively.}
\includegraphics[scale=0.8]{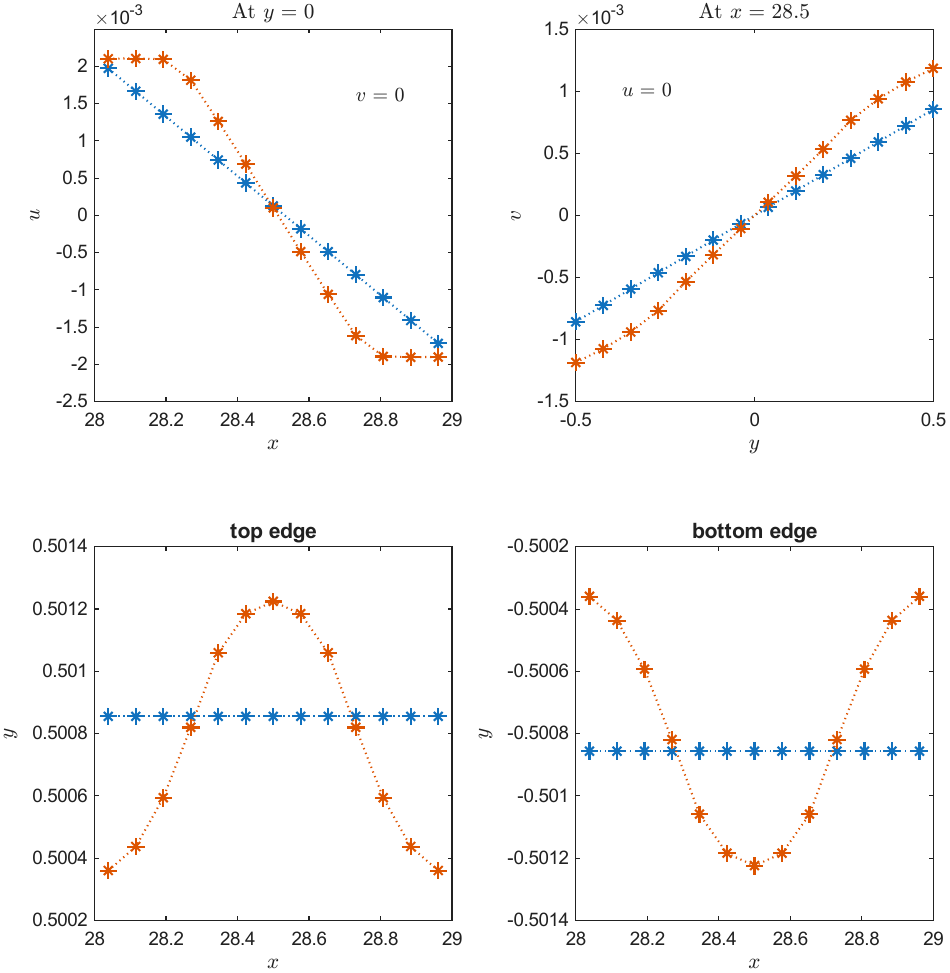}  
\end{figure}

\cref{fig:uvtbedgesEs1001}(blue) shows computational results that demonstrate this linear dependence in the centre patch of the domain~(\(28\le x\le 29\), all quantities non-dimensional).
As the beam is compressed in the axial direction, we observe linear compression in~\(u\), but under this compression the beam thickens in the transverse direction, leading to a linear change in~\(v\) across the cross-section of the beam~(\cref{fig:uvtbedgesEs1001}, top row). 
The ratio between the slope of the~\(v\)-profile and the slope of the~\(u\)-profile in the uniform stress region is approximately~\(0.3\), which is the prescribed Poisson's ratio for the material.  
\cref{fig:uvtbedgesEs1001}(blue, bottom row) shows that the top and bottom edges of the beam remain straight when under this axial strain, indicating \text{an absence of buckling.}

\begin{figure}
\centering
\caption{\label{fig:SI001Stress}
Axial normal stress~\(\sigma^{xx}\) (top row), and shear stress~\(\sigma^{xy}\) (bottom row), in the three leftmost patches of the soft-inclusion beam with~\(E_s=0.01\) under axial compression strain~\(\varepsilon^e=0.4\%\) (so that \(0.11<x<56.89\)), obtained from the full-domain computation~(\(r=1\)). }
\begin{tabular}{@{}c@{}c@{}c@{}}
\includegraphics{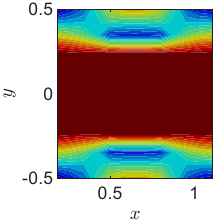} & 
\includegraphics{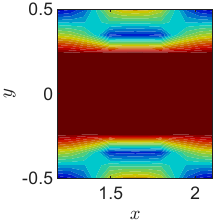} & 
\includegraphics{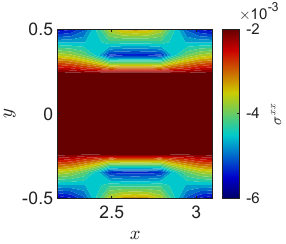}\\
 \includegraphics{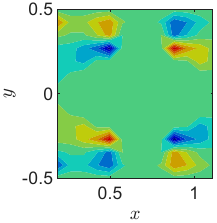} & 
 \includegraphics{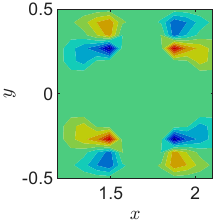} & 
 \includegraphics{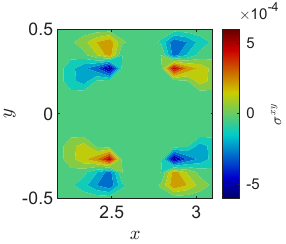}\\
\end{tabular}
\end{figure}%

For axial strain~\(\varepsilon^e=0.4\%\)  in a soft-inclusion beam, \cref{fig:SI001Stress} shows the axial normal stress~\(\sigma^{xx}\) and the shear stress~\(\sigma^{xy}\) in the three leftmost patches. 
In this case, the difference in stiffness between the matrix material~(\(E=1\)) and the soft inclusions~(\(E_s<1\)) generates a non-uniform stress field that is spatially periodic (corresponding to the period of the beam structure), except in the boundary layers close to the ends of the beam. 
Near the axis of the beam~(\(y\approx 0\)), the normal stress~\(\sigma^{xx}\) is significantly smaller in magnitude than  near the top and bottom edges of the beam, due to the inclusions being less stiff than the surrounding matrix. 
In contrast to the homogeneous beam, the shear stress~\(\sigma^{xy}\) is no longer zero and is relatively large in regions about the soft inclusions. 
This shear effect is induced by the abrupt spatial change in elasticity across the inclusion-matrix interface, but is only significant inside the stiffer matrix material.  
The resulting stress localisation is consistent with classical elasticity solutions for circular elastic inclusions, which show that a mismatch in elastic properties between an inclusion and its surrounding matrix induces local variations in the stress field~\citep{Goodier1933}. 
Similar behaviour has been reported for weak circular inclusions subjected to uniaxial compression, where pronounced stress concentrations develop around the inclusion--matrix interface and depend strongly on the stiffness contrast between the two materials~\citep{WuWong2013, Zhu2022}.
The non-uniform stress field gives rise to \text{nonlinear displacement fields.}

\Cref{fig:uvtbedgesEs1001}(red) illustrates the nonlinear displacements in the centre patch \((28\le x\le 29)\) of the soft-inclusion beam.
Within both the matrix material and the soft inclusion, the horizontal displacement field~\(u\) is approximately linear, although with distinctly different slopes in the two materials. 
The steeper slope of the horizontal displacement within the inclusion is due the normal strain~\(\varepsilon^{xx}\) being larger inside this soft inclusion compared to the surrounding matrix material. 
Across the inclusion-matrix interface the  horizontal displacement has a nonlinear transition, resulting in a variable displacement gradient~\(\partial u/\partial x\) and variations in local strains and stresses.
Under axial compression the beam thickens in the transverse direction, with the largest increases in the vertical displacement field~\(v\) about the centre of the soft inclusion.
Consequently, the top and bottom beam boundaries are rippled, with maximum displacements aligned with the centres of the inclusions. 
The top and bottom boundaries are symmetric \text{about the beam axis.}  

\Cref{fig:uvtbedgesEs1001} shows displacements from both full-domain computations (dotted lines) and from computations with \patRats~\cref{Epatrat} of~\(r=0.5\) and~\(0.25\) ($+$~and~$\times$). 
The patch computations are nearly indistinguishable from the full domain computations, whether considering  the homogeneous beam~(\(E_s=1\)) or the soft-inclusion beam~(\(E_s=0.01\)).
These results demonstrate that the patch scheme accurately captures the unbuckled equilibrium states, despite a  substantial reduction in the computational domain, and validating theoretical results of~\cite{Bunder2020a}.
To explore accuracy in  less trivial predictions of the patch scheme, \cref{subSec:Emodes} examines some dynamical and bifurcation properties of the eigenvalues and eigenvectors of \text{the multiscale beam system.}

\subsection{Onset of buckling}
\label{subSec:Emodes}

Under increasing axial compression strain, the unbuckled configuration remains an equilibrium state of the beam owing to the up--down symmetry of the beam geometry and elasticity. 
Nevertheless, this unbuckled equilibrium undergoes a transition from stable to unstable at some critical axial compression, characteristic of the pitchfork bifurcation associated with the buckling of symmetric beams~\citep{GhayeshFarokhi2015}. 
According to bifurcation theory, the stability transition occurs when one or more eigenvalues of the Jacobian matrix~\cref{eq:JabMatrix} acquire positive real parts. 
Beyond this critical compression, even infinitesimal perturbations can trigger a rapid evolution away from the now unstable unbuckled configuration towards a \text{new buckled equilibrium.}

For a two dimensional beam, two distinct types of deformation modes generally exist: compressive modes and bending modes. 
The former correspond to deformations predominantly along the beam axis, whereas the latter are associated with transverse deflections and bending of the beam and thus are associated with buckling. 
Here it is the long-wavelength bending modes that are of interest since these are the first to destabilise as \text{the axial compression increases.}

\paragraph{Homogenous beams}
\begin{figure}
\centering
\caption{\label{fig:EVal_1}
For a homogeneous beam~(\(E_s=1\)), plot the variation in patch-scheme eigenvalues as a function of the axial compression strain, \(\varepsilon^e\leq0.4\%\). 
The chosen \patRats\ are~\cref{Epatrat}: (a)~\(r=1\); (b)~\(r=0.5\); and (c)~\(r=0.25\).
}
\begin{enumerate}[nosep,label=(\alph*)]
\item Full-domain (\(r=1\))   \\
\includegraphics[scale=0.8]{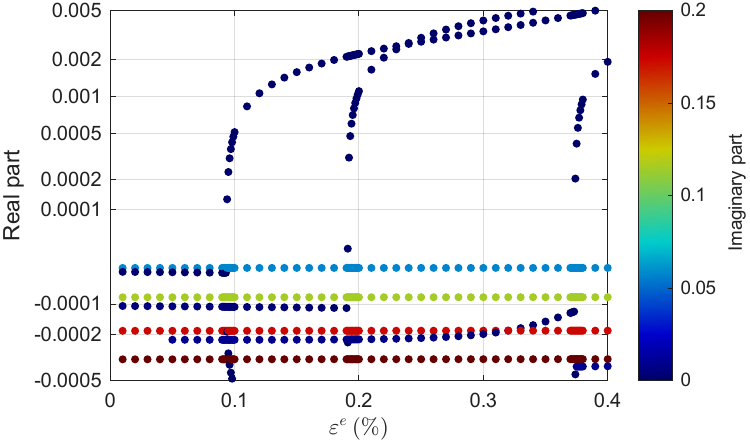}    
\item Half-domain (\(r=0.5\)) \\
\includegraphics[scale=0.8]{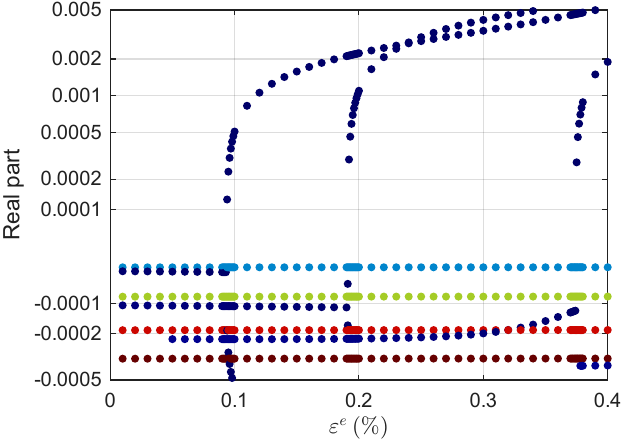}    
\item Quarter-domain (\(r=0.25\)) \\
\includegraphics[scale=0.8]{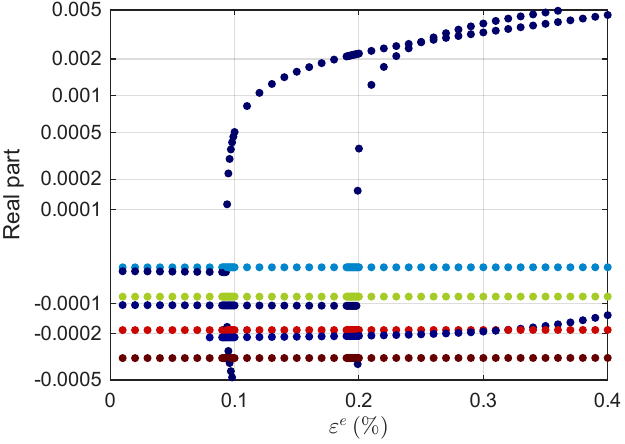}
\end{enumerate}
\end{figure}%
For the homogeneous beam with axial compression strains up to~\(\varepsilon^e=0.4\%\), \Cref{fig:EVal_1} plots those eigenvalues of the Jacobian matrix~\cref{eq:JabMatrix} with small real~part, as computed using different \patRats~\cref{Epatrat}.  
Since the eigenvalues are generally complex, for simplicity, the vertical coordinate represents the real-parts and the colour indicates the imaginary-parts.  
The real-part characterises the stability state of the corresponding eigen-mode, whereas the imaginary-part represents the frequency of the mode.  
Only eigenvalues with real-part close to zero are displayed, since buckling instability occurs when the real-part transitions from \text{negative to positive.}

\Cref{fig:EVal_1} plots both bending and compressive eigenvalues.
The long wavelength bending eigen-modes have small imaginary-parts (low frequencies) and so are shown by dark blue dots.  
The compressive eigen-modes have larger imaginary-parts (higher frequencies) and are indicated by green--red dots.  
The bending eigenvalues have almost constant negative real-parts for small axial strains, but undergo dramatic increases to positive real-parts at threshold applied stains--thresholds which differ for different modes.  
In contrast, the compressive eigenvalues have almost constant real parts for all \text{computed axial strains.}

\begin{figure}
\centering
\caption{The first three bending eigen-modes of the homogeneous beam~(\(E_s=1\)) under controlled axial strain~\(\varepsilon^e=0.02\%\).
For each case there are two plots corresponding to: (top) half-domain patches~\(r=0.5\); and (bottom) quarter-domain patches~\(r=0.25\).  
The colour indicates the normal stress~\(\sigma^{xx}\).}\label{fig:EShape_1}
\begin{enumerate}[nosep,label=(\alph*)]
\item Unimode  
\item[]\includegraphics[scale=0.8]{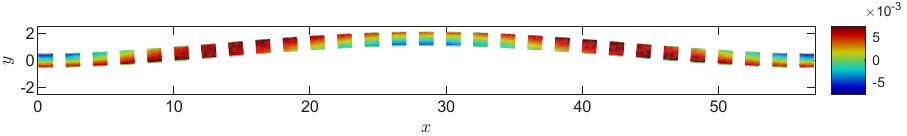}    
\item[]\includegraphics[scale=0.8]{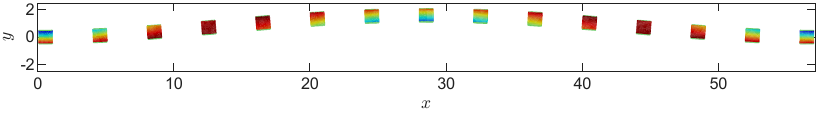}    
\item Bimode 
\item[]\includegraphics[scale=0.8]{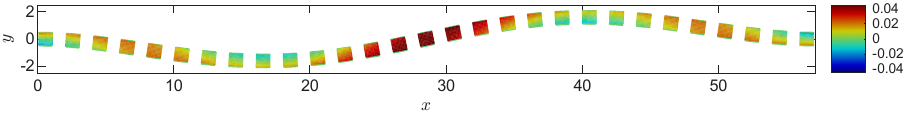}    
\item[]\includegraphics[scale=0.8]{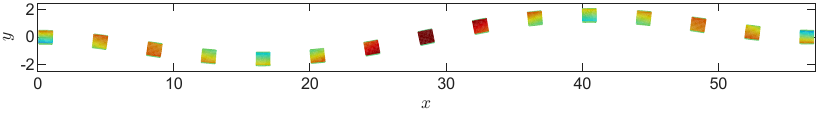}    
\item Trimode 
\item[]\includegraphics[scale=0.8]{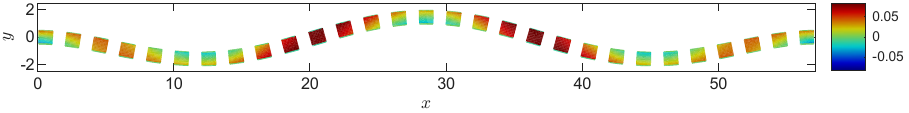}    
\item[]\includegraphics[scale=0.8]{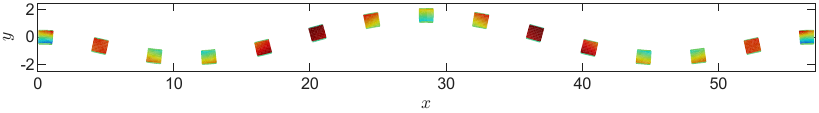}    
\end{enumerate}
\end{figure}%

In the full-domain computation~(\(r=1\)), the first three bending modes to undergo stability changes do so  when the axial strain~\(\varepsilon^e\) exceeds~\(0.093\%\), \(0.191\%\), and~\(0.373\%\)~(\Cref{fig:EVal_1}(a)).
These critical strain thresholds are hereafter denoted by~\(\varepsilon_{\Cr}\).   
\cref{fig:EShape_1} illustrates the corresponding beam modes, computed using~\(r=0.5\) and~\(r=0.25\), and shows that they are the three lowest frequency bending modes.
We name these three modes the unimode, bimode, and trimode, based on the number of buckles in the beam.
The three observed instability transitions occur in order of lowest to highest frequency (or  number of buckles) as~\(\varepsilon^e\) increases.
The critical axial normal stresses associated with these transitions are \(\sigma_{\Cr}^{xx}=1.03\times10^{-3}\), \(2.10\times10^{-3}\), \text{and~\(4.11\times10^{-3}\), respectively.}  

The critical axial compression strain thresholds~\(\varepsilon_{\Cr}\) and associated axial normal stresses~\(\sigma_{\Cr}\) for reduced domain~\(r<1\) computations have small differences compared to the full-domain computation~(\(r=1\)) that increase as~\(r\) decreases~(\Cref{fig:EVal_1}(b)--(c)). 
These differences are due to the naturally increasing \Ord{H^m}~errors (\cref{subSec:PatchCoup}) because the spacing~\(H\) between patches increases with decreasing~\(r\).
\Cref{tab:StrainThresh} lists these critical values for the homogeneous beam~(\(E_s=1\)) and \patRats\ \(r=1,0.5,0.25\).
Typically, reducing~\(r\), thus increasing~\(H\), overestimates the critical strain thresholds.
For the half-domain computations~\(r=0.5\), the overestimation is quite small, with critical strains within three decimal percent-places of those computed on the full-domain.
For the quarter-domain~\(r=0.25\), the critical strain differences are a little larger. 
In general, the patch spacing has less influence on the important unimode stability threshold than it does \text{on the higher modes.}  

\begin{table}
    \centering
    \caption{\label{tab:StrainThresh}
    Critical strain thresholds~\(\varepsilon_{\Cr}\) and corresponding critical stresses~\(\sigma_{\Cr}\) for the  first three bending eigen-modes of the homogenous beam~\(E_s=1\), and the soft-inclusion beams~\(E_s=0.001\)--\(0.1\), computed with \patRats~\cref{Epatrat} \(r=1,0.5,0.25\). }\
\begin{equation*}
    \begin{array}{|cc|cc|cc|cc|}
     \hline
           \multirow{2}{*}{\(E_s\)}& \multirow{2}{*}{\text{Modes}} &\multicolumn{2}{c|}{r=0.25}      & \multicolumn{2}{c|}{r=0.5}     & \multicolumn{2}{c|}{r=1} \\
              & & \varepsilon_\Cr  & \sigma_\Cr \cdot10^{3}& \varepsilon_\Cr  & \sigma_\Cr \cdot10^{3}& \varepsilon_\Cr  & \sigma_\Cr \cdot10^{3} \\
      \hline
       \multirow{3}{*}{1}                 & \text{1st}   & 0.093\%  & 1.024 & 0.093\% & 1.023 & 0.093\% & 1.025 \\
                                                    & \text{2nd} & 0.199\%  & 2.190 & 0.191\% & 2.100 & 0.191\% & 2.103 \\
                                                    & \text{3rd}  & 0.441\%  & 4.847 & 0.375\% & 4.118 & 0.373\% & 4.105 \\
      \hline
       \multirow{3}{*}{\(10^{-1}\)}  & \text{1st}   & 0.080\% & 0.544 & 0.081\% & 0.569 & 0.081\% & 0.568 \\
                                                    & \text{2nd} & 0.169\% & 1.149 & 0.165\% & 1.159 & 0.163\% & 1.143 \\
                                                    & \text{3rd}  & 0.362\% & 2.464 & 0.335\% & 2.351 & 0.323\% & 2.313 \\
     \hline
      \multirow{3}{*}{\(10^{-2}\)}  & \text{1st}   & 0.078\% & 0.468 & 0.079\% & 0.493 & 0.079\% & 0.493 \\
                                                   & \text{2nd} & 0.162\% & 0.987 & 0.157\% & 0.992 & 0.157\% & 0.980 \\
                                                   & \text{3rd}  & 0.345\% & 2.067 & 0.306\% & 2.008 & 0.303\% & 1.933 \\    
    \hline
     \multirow{3}{*}{\(10^{-3}\)}   & \text{1st}   & 0.077\% & 0.459 & 0.077\% & 0.459& 0.077\% & 0.474\\
                                                   & \text{2nd} & 0.163\% & 0.967 & 0.156\% & 0.966 & 0.155\% & 0.954 \\
                                                   & \text{3rd}  & 0.343\% & 2.028 & 0.306\% & 2.010 & 0.304\% & 1.906\\
     \hline      
\end{array}\end{equation*}
\end{table}

\paragraph{Heterogeneous soft-inclusion beams}
\begin{figure}
\centering
\caption{\label{fig:EVal_001}
For the soft-inclusion beam with~\(E_s=0.01\), plot the variation in eigenvalues versus axial compression~\(\varepsilon^e\leq0.4\%\). 
The chosen \patRats~\cref{Epatrat} are: (a)~\(r=1\); (b)~\(r=0.5\); and (c)~\(r=0.25\).}
\begin{enumerate}[nosep,label=(\alph*)]
\item Full-domain (\(r=1\)) 
\item[]\includegraphics[scale=0.8]{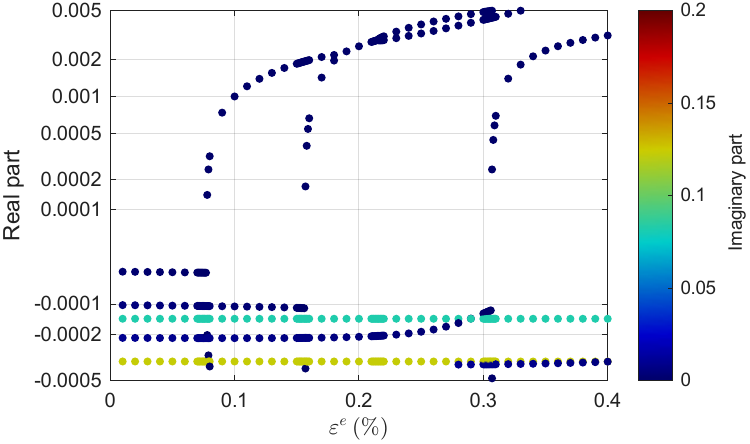}    
\item Half-domain (\(r=0.5\))  
\item[]\includegraphics[scale=0.8]{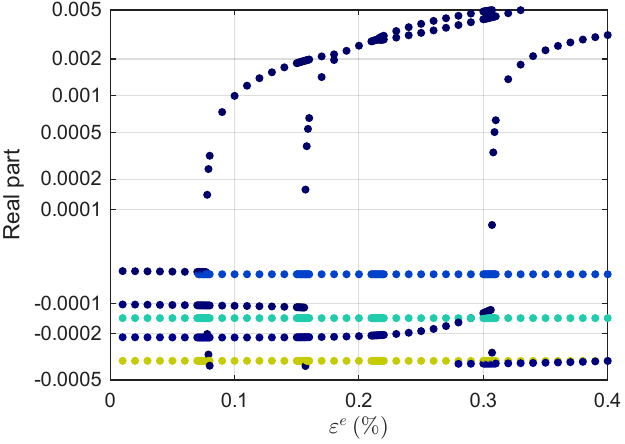}    
\item Quarter-domain (\(r=0.25\)) 
\item[]\includegraphics[scale=0.8]{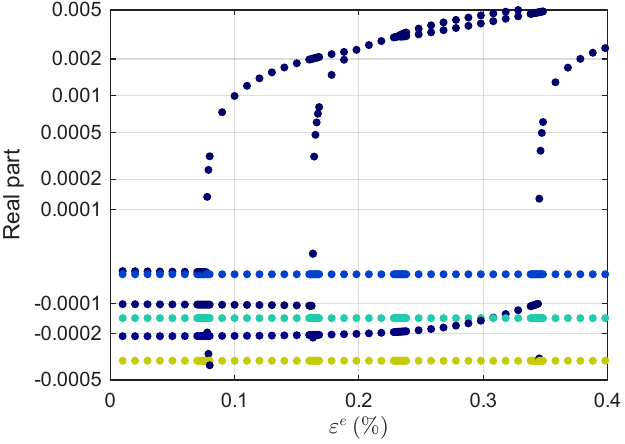}
\end{enumerate}
\end{figure}%

\begin{figure}
\centering
\caption{\label{fig:EShape_001}
The first three bending eigen-modes of the soft-inclusion beam with~\(E_s=0.01\) under compression with~\(\varepsilon^e=0.02\%\).
For each case there are two plots corresponding to: (top) half-domain patches~\(r=0.5\); and (bottom) quarter-domain patches~\(r=0.25\).  
The colour indicates the normal stress~\(\sigma^{xx}\).}
\begin{enumerate}[nosep,label=(\alph*)]
\item Unimode  
\item[]\includegraphics[scale=0.8]{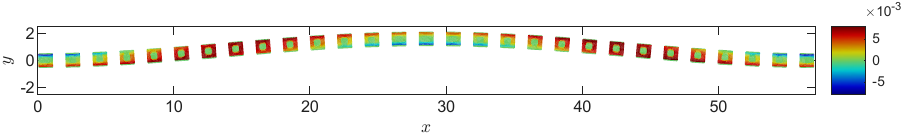}    
\item[]\includegraphics[scale=0.8]{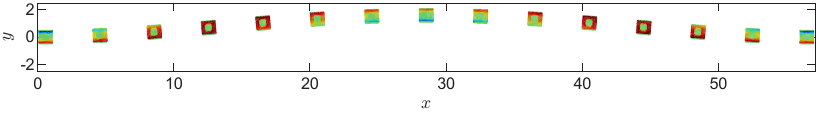}    
\item Bimode
\item[]\includegraphics[scale=0.8]{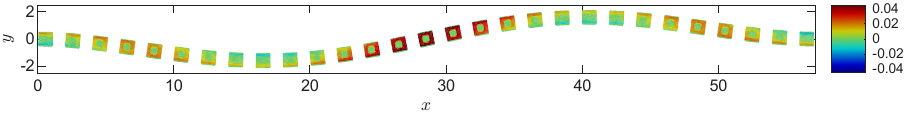}    
\item[]\includegraphics[scale=0.8]{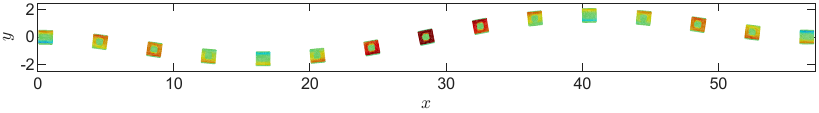}    
\item Trimode
\item[]\includegraphics[scale=0.8]{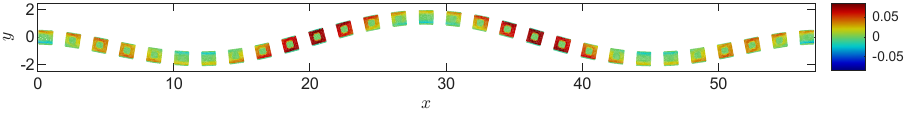}    
\item[]\includegraphics[scale=0.8]{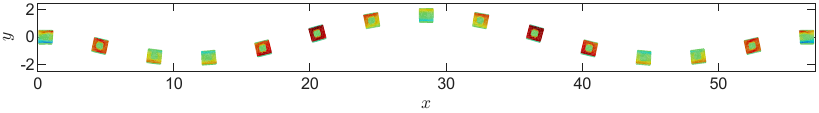}    
\end{enumerate}
\end{figure}

\Cref{fig:EVal_001}  illustrates the small real part eigenvalues of the Jacobian matrix for the soft-inclusion beam with Young's modulus \(E_s=0.01\) under axial compression, computed using the \patRats~\cref{Epatrat} \(r=1\), \(0.5\), and~\(0.25\).  
The characteristics of the eigenvalues are qualitatively similar to those observed for the homogeneous beam~(\cref{fig:EVal_1}).  
Specifically, the compressive modes remain stable throughout the compression process, whereas the three lowest-frequency bending modes are observed to undergo stability transitions once the axial strain exceeds critical thresholds.  
Compared with the homogeneous beam, the soft-inclusion beam undergoes stability transitions at substantially smaller axial strains, and critical stresses are approximately one-half of those obtained for the homogeneous beam.  
\cref{fig:EShape_001} displays the corresponding eigen-modes of the soft-inclusion beam with colours indicating normal stress~\(\sigma^{xx}\).  
Although the global deformation patterns are similar to those of the homogeneous beam, the presence of soft inclusions introduces pronounced localised internal structure, as shown by the small-scale variations in normal stress, particularly \text{near interface boundaries.} 

\Cref{tab:StrainThresh} summarises the critical strain thresholds~\(\varepsilon_{\Cr}\), and the corresponding critical stresses~\(\sigma_{\Cr}\) for  the soft-inclusion beams with~\(E_s=0.1\)--\(0.001\).  
Decreasing the Young's modulus of the soft inclusions reduces the critical strain threshold of the bending modes, as is expected, given that a smaller Young's modulus implies a greater propensity for deformation.  
The soft-inclusion beams also exhibit lower critical stresses compared to the homogenous beam, demonstrating that the inclusions reduce the beam's overall structural stiffness.
The influence of the patch spacing~\(H\) on the stability transitions of the soft-inclusion beams is qualitatively similar to that observed for the homogeneous beam.  
In particular, larger values of spacing~\(H\) (smaller \patRat~\(r\)) typically lead to larger overestimations of critical strain thresholds of the bending modes, with this difference more pronounced for higher frequency modes, although the overestimations are \text{not particularly large.}

\subsection{Buckled beams}
\label{SSbb}

In a compressed beam, buckling occurs through a symmetry-breaking instability in which the straight equilibrium loses stability and the beam evolves towards one of the two transversely deflected configurations~\citep{Wan2020}.
The elastic strain energy stored under axial compression is released as kinetic energy and drives the beam away from its straight equilibrium configuration, leading to a transverse deflection of the beam and the formation of a buckled shape.
This section explores predictions of equilibrium buckled states and the accuracy of \text{the patch scheme.}

To quantify the degree of buckling, the transverse deflection~\(\zeta(x)\) is defined as the vertical displacement of the beam axis from~\(y=0\). 
For an unbuckled beam~\(\zeta(x)=0\) for~\(0\le x\le L\).
Following the onset of buckling, the maximum transverse deflection~\(\zeta_{\max}\) is used to characterise the buckled state.
Owing to the left-right symmetry of both the beam and its axial compression, \(\zeta_{\max}\)~typically occurs half-way along the beam at~\(x=L/2\) \text{(the unimode case).}
     
In numerical computations, to determine the buckled equilibrium, a small initial disturbance is introduced to break the symmetry of the system. 
We use the first two bending modes (uni- and bimode, \cref{fig:EShape_1,,fig:EShape_001}) to seed the disturbance. 
Accordingly, instead of setting the initial displacement to~\(\uv_0=0\), the displacement field is initialised as~\(\uv_0 = \alpha\uv_{k}\),
where~\(\uv_{k}\) denotes the displacement field associated with the \(k\)th bending mode, and~\(\alpha\) is a small prescribed disturbance magnitude.
Once a buckled equilibrium is found, the axial strain~\(\varepsilon^e\) is incremented by small amounts, and the nearby buckled equilibrium is sought for each \text{increment in~\(\varepsilon^e\).}

\begin{figure}
	\centering
	\caption{\label{fig:Buckled_1}
	Buckled configurations of the homogeneous beam with~\(E_s=1\) and maximumn transverse deflection: (a)~\(\zeta_{\max}=0.4\); and (b)~\(\zeta_{\max}=0.9\).
		For each case there are two plots corresponding to: (top)~\(r=0.5\); and (bottom)~\(r=0.25\). 
		The colour indicates the axial normal stress \(\sigma^{xx}\).}
	\begin{enumerate}[nosep,label=(\alph*)]
		\item \(\zeta_{\max}=0.4\)
		\item[]\includegraphics[scale=0.8]{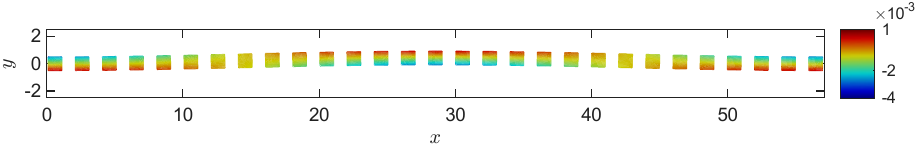} 
		\item[]\includegraphics[scale=0.8]{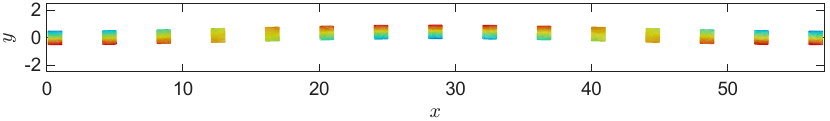} 
		\item \(\zeta_{\max}=0.9\)
		\item[]\includegraphics[scale=0.8]{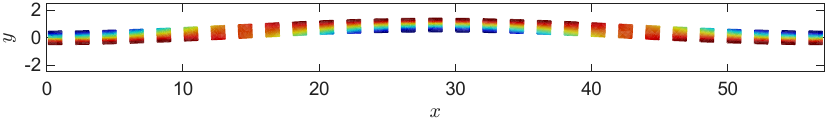} 
		\item[]\includegraphics[scale=0.8]{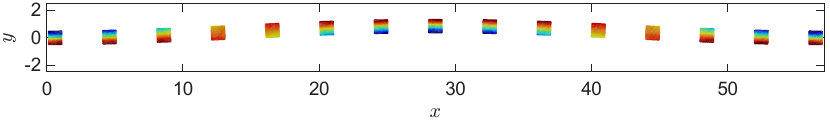} 
	\end{enumerate}
\end{figure}

\paragraph{Homogeneous beam} 
\Cref{fig:Buckled_1} shows the beam equilibria for buckled states with~\(\zeta_{\max}=0.4\) and~\(0.9\), as predicted from the patch scheme with \patRats~\cref{Epatrat} \(r=0.5\) and~\(0.25\).
For the two \patRats, the buckled beam shapes are almost the same, and the stress distributions within the patches are in close agreement (as indicated by the colours).
The primary difference in the buckling for the different \patRats\ is that a specific~\(\zeta_{\max}\) is achieve at slightly different axial strains.
For \patRat\ \(r=0.5\), the deflections~\(\zeta_{\max}=0.4\) and~\(0.9\) occur at axial strains~\(\varepsilon^e=0.103\%\) and~\(0.128\%\), respectively;
 whereas for  \(r=0.25\), these deflections occur at axial strains~\(0.104\%\) and~\(0.134\%\), respectively.
Thus, as~\(r\) decreases, slightly larger axial strains are typically predicted to produce a specific axial deflection, which is consistent with \Cref{tab:StrainThresh} where buckling instability computed with smaller~\(r\) is usually predicted to occur at \text{larger critical axial strain.} 

\begin{figure}
\centering
\caption{\label{fig:VMaxCrStress}
Maximum transverse deflection~\(\zeta_{\max}\), and critical buckling stress~\(\sigma_\Cr\), for compressed soft-inclusion beams with: (a)~\(E_s=1\) (homogeneous); (b)~\(E_s=0.1\); and (c)~\(E_s=0.01\).
	Each case is computed with three \patRats: full domain~\(r=1\) (blue dotted curves); half-domain~\(r=0.5\) (red~\(+\)); and quarter-domain~\(r=0.25\) (orange~\(\times\)).}
\begin{enumerate}[nosep,label=(\alph*)]
\item \(E_s=1\)
\item[]\includegraphics[scale=0.8]{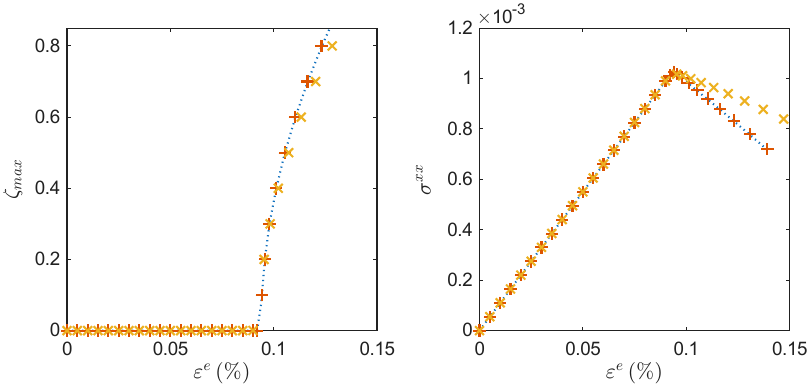} 
\item \(E_s=0.1\)
\item[]\includegraphics[scale=0.8]{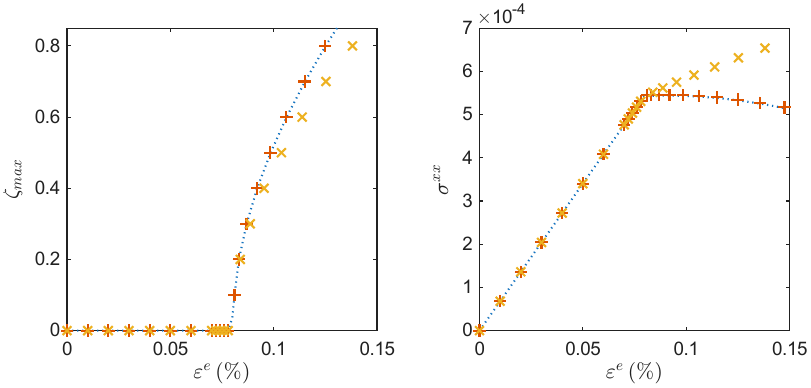} 
\item \(E_s=0.01\)
\item[]\includegraphics[scale=0.8]{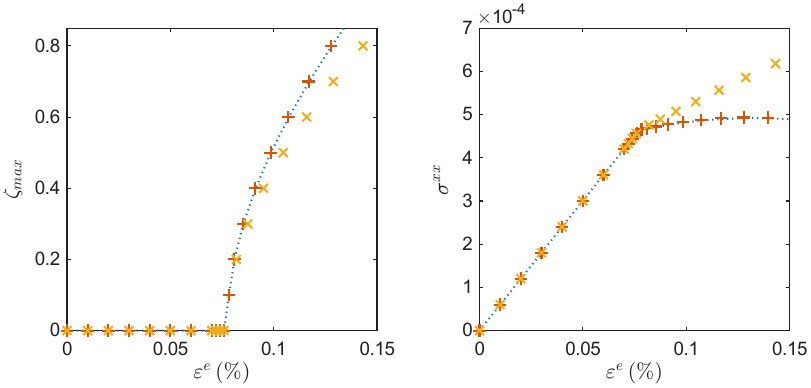}  
\end{enumerate}
\end{figure}%

\Cref{fig:VMaxCrStress}(a) plots the evolution of the maximum deflection~\(\zeta_{\max}\) and the axial normal stress \(\sigma^{xx}\) due to axial compression strain~\(\varepsilon^e\). 
The deflection deviates from zero once the axial strain exceeds approximately~\(\varepsilon_{cr}=0.093\%\) (the first critical axial strain threshold in~\Cref{tab:StrainThresh}). 
Prior to the onset of buckling, the normal stress~\(\sigma^{xx}\) increases linearly with the applied axial compressive strain~\(\varepsilon^e\). 
Once buckled, at approximately~\(\varepsilon_{cr}=0.093\%\) and~\(\sigma_{cr}=1.025\times10^{-3}\) (\Cref{tab:StrainThresh}), \(\sigma^{xx}\) begins to decrease as the compression continues to increase.
The deflection and stress curves predicted by half-domain patches, \(r=0.5\), closely follow those of the full-domain computation with excellent accuracy. 
For the quarter-domain patches, \(r=0.25\), the deflection and stress curves are consistent with the full-domain results until buckling, after which there is an underestimation in the deflection and overestimated in the stress for a given axial strain, again consistent with \Cref{tab:StrainThresh}.
The overestimation of the stress in the~\(r=0.25\) buckled beam is understandable, as the increased patch spacing~\(H\) limits the ability of the patch scheme to accurately resolve the increasingly nonlinear beam deformations that develop in the post-buckling regime.
Moreover, an overestimation of the stress indicates an underestimation of the released elastic energy as kinetic energy, and thus leads to an underestimation of \text{the beam deflection.}

The critical stess~\(\sigma_{\Cr}=1.025\times 10^{-3}\) (\Cref{tab:StrainThresh} and \Cref{fig:VMaxCrStress}(a)) at which buckling first occurs is associated with a critical force or load~\(F^e_{\Cr}\).
From classical Euler buckling theory \cite[e.g.,][]{Euler1947, Avcar2014}, the critical load at which a straight, slender, homogeneous beam loses stability is
\begin{equation}
F^e_{\Cr} := {\pi^2 E I}/{(K L)^2},\label{eq:EulerCriBL}
\end{equation}
where~\(E\) is the Young's modulus, \(I\) is the second moment of area of the beam cross-section, and~\(L\) is the beam length. 
The parameter~\(K\) is the effective length factor accounting for the boundary conditions:  \(K=1\) for pinned-pinned beams; \(K=0.5\) for fixed-fixed beams; and~\(K=2\) for fixed-free beams. 
In the present study, both beam ends are fixed, and hence here~\(K=0.5\). 
Substituting the non-dimensional quantities~\(E=1\), \(I=1/12\), corresponding to a unit beam height, and length~\(L=57\) into~\cref{eq:EulerCriBL}, and dividing by the cross-sectional area \(A=1\), yields the corresponding critical stress for buckling: \(\sigma_{\Cr} = {F^e_{\Cr}}/{A}= 1.025\times 10^{-3} \),
which agrees to three decimal places with the value obtained from the full-domain \text{and patch computations. }

\paragraph{Heterogeneous soft-inclusion beam}

\begin{figure}
\centering
\caption{Buckled configurations of the soft-inclusion beam with \(E_s=0.01\) at maximumn transverse deflection: (a)~\(\zeta_{\max}=0.4\); and (b)~\(\zeta_{\max}=0.9\). 
	For each case there are two plots corresponding to: (top) \(r=0.5\); and (bottom) \(r=0.25\). 
	The colour indicates the axial normal stress \(\sigma^{xx}\).}\label{fig:Buckled_001}
\begin{enumerate}[nosep,label=(\alph*)]
\item \(\zeta_{\max}=0.4\)
\item[]\includegraphics[scale=0.8]{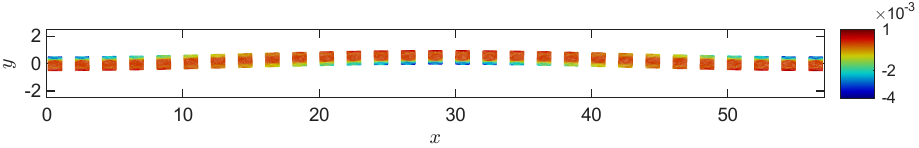} 
\item[]\includegraphics[scale=0.8]{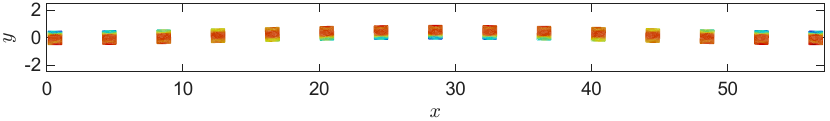} 
\item \(\zeta_{\max}=0.9\)
\item[]\includegraphics[scale=0.8]{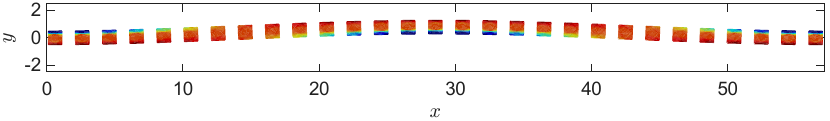} 
\item[]\includegraphics[scale=0.8]{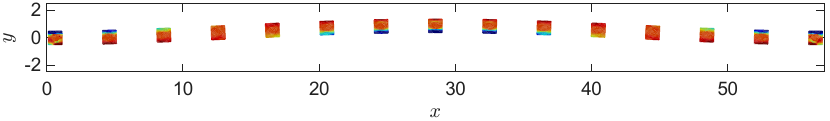} 
\end{enumerate}
\end{figure} 

\begin{figure}
\centering
\caption{\label{fig:CrStrainStress}
Critical strain~\(\varepsilon_{\Cr}\) and critical stress~\(\sigma_{\Cr}\) for buckling instability in the compressed soft-inclusion beams predicted with \patRats\ \(r=1,0.5,0.25\). }
\includegraphics[scale=0.9]{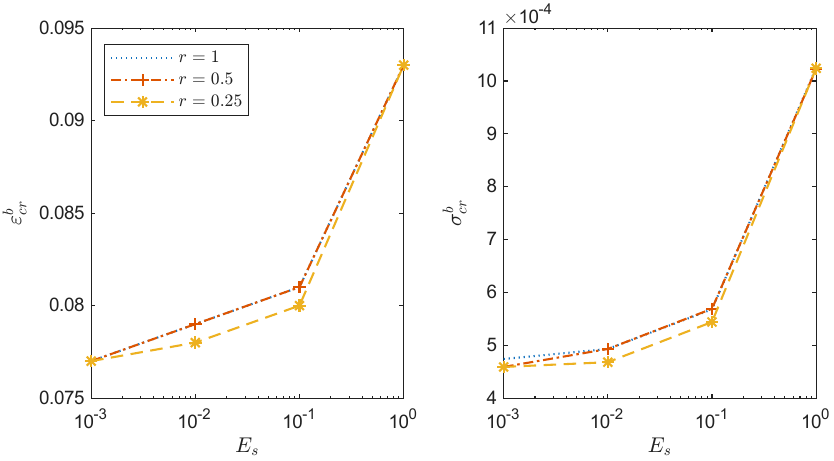}    
\end{figure}

\Cref{fig:VMaxCrStress}(b)--(c) plots the variation of the maximum deflection~\(\zeta_{\max}\) and the axial normal stress~\(\sigma^{xx}\) due to axial compression strain~\(\varepsilon^e\) for the soft-inclusion beams with~\(E_s=0.01,0.1\).
As the stiffness of the soft inclusions decrease with Young's modulus~\(E_s\), the beam undergoes buckling instability at lower critical strains and stresses than the homogeneous beam (also shown in~\Cref{tab:StrainThresh}).
The influence of the soft inclusions is evident not only in the critical buckling values, but also in the post-buckling response. 
After bucking, transverse deflections increase more slowly and stresses decreases less rapidly with decreasing~\(E_s\), as strain tends to be more local to the inclusions, rather than a macroscale property of the beam.
The soft inclusions accommodate a greater proportion of the imposed compression and associated elastic energy through local deformation, thereby mitigating the global post-buckling response.
Consequently, beams with softer inclusions exhibit a milder post-buckling behaviour, characterised by smaller transverse deflections and \text{weaker stress relaxation.}
 
\Cref{fig:CrStrainStress} plots the critical strain and critical stress for the buckling instability over a wide range of soft-inclusion stiffness values, \(E_s=0.001\)--\(1\).
Both critical quantities decrease rapidly as the stiffness of the soft inclusions is reduced from~\(E_s=1\) to~\(E_s=0.1\).
As the stiffness is further decreased, the rate of reduction becomes progressively smaller.
Consequently, the critical curves exhibit an asymptotic trend as~\(E_s\) approaches zero.
This asymptotic behaviour is expected because the limit~\(E_s\to0\) represents a porous structure with void-like inclusions. 
Previous studies of porous beams similarly demonstrate a reduction in
critical buckling resistance associated with increasing porosity and the
corresponding loss of effective stiffness~\citep{Chen2015, Wu2024}.
In the present beam, as~\(E_s\) approaches zeor, the inclusions carry progressively less stress and approach the limiting behaviour of voids; consequently, further reductions in~\(E_s\) have a diminishing influence on \text{the global buckling response.} 

\paragraph{Bifurcation digram}
Buckling is the sudden transition from one stable equilibrium configuration of a beam to another, driven by an increase in the axial compression strain~\(\varepsilon^e\) beyond a critical value. 
Mathematically, buckling usually corresponds to a pitchfork bifurcation when some bifurcation parameter, here the applied strain~\(\varepsilon^e\), is increased beyond its critical value, \text{here denoted~\(\varepsilon_{\Cr}\).}

\begin{figure}
\centering
\caption{\label{fig:BiDiagram}
Bifurcation diagrams for compressed soft-inclusion beams computed predicted by the patch scheme with two \patRats: \(r=0.5\) (blue) and~\(r=0.25\) (red).
The green branch represents the unbuckled (straight) equilibrium solution. 
Solid and dashed lines respectively denote stable and \text{unstable solution branches.} }
\begin{enumerate}[nosep,label=(\alph*)]
\item \(E_s=1\)
\item[] \includegraphics[scale=0.8]{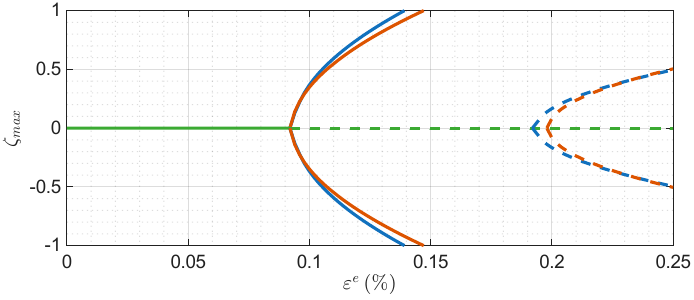} 
\item \(E_s=0.1\)
\item[] \includegraphics[scale=0.8]{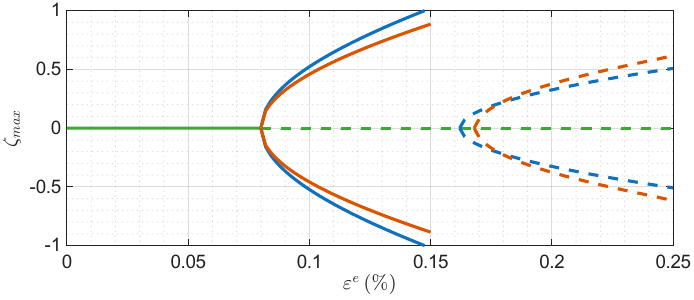} 
\item \(E_s=0.01\)
\item[] \includegraphics[scale=0.8]{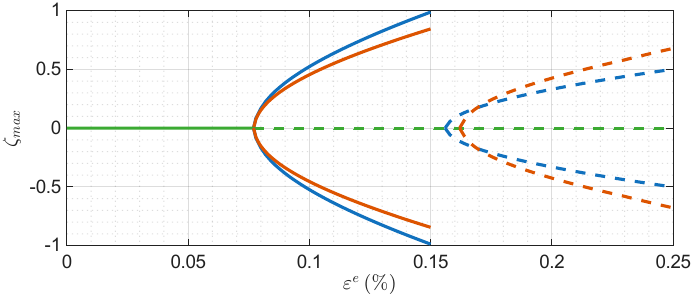} 
\item \(E_s=0.001\)
\item[] \includegraphics[scale=0.8]{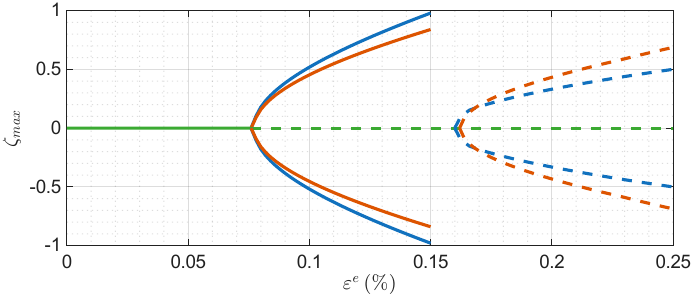} 
\end{enumerate}
\end{figure}

\Cref{fig:BiDiagram} presents bifurcation diagrams for compressed soft-inclusion beams with different Young's modulus~\(E_s\).
Solid and dashed lines denote stable and unstable equilibrium branches, respectively.
The blue and red curves correspond to the \patRats~\cref{Epatrat} \(r=0.5\) and~\(0.25\), respectively.
A new pitchfork bifurcation appears with each new mode of buckling, although it is only the unimodal buckling branches that are stable for~\(\varepsilon^e>\varepsilon_{\Cr}\), with bimodal and higher eigen-modes always unstable.
\Cref{fig:BiDiagram} shows that as~\(E_s\) decreases,  unimodal and bimodal buckling branches usually occur at lower axial strains, although the difference is small for the smallest~\(E_s\), \text{in agreement with~\cref{fig:CrStrainStress}.}

The bifurcation diagrams provide valuable insight into the buckling behaviour of the beams, including the critical buckling points, the number of equilibrium branches, and the evolution of the post-buckling deflection.
They can be obtained either from direct buckling simulations or by numerical continuation methods, which systematically trace equilibrium branches while determining their stability through analysis of the Jacobian matrix.
Several software packages provide numerical continuation and bifurcation-tracking capabilities, including Matlab package \textsc{matcont}, for systems of small \text{or intermediate size.}
    
\subsection{Order of the interpolation}
\label{SSoi}

An important parameter in the patch scheme is the order~\(m\) of the interpolating polynomial~(\cref{subSec:PatchCoup}), which directly influences the accuracy of the patch coupling as its error is~\Ord{H^m}.
Consequently, the order~\(m\) influences the overall accuracy of the multiscale simulation. 
In practice, low-order interpolations~(\(m=2\) or~\(4\)) are often sufficient for linear systems or systems exhibiting gentle spatial variations, whereas higher-order interpolations (\(m\ge 6\)) are generally more effective for nonlinear systems or systems with \text{strong spatial heterogeneity.}

To demonstrate the influence of the interpolation order on the accuracy of the patch scheme, particularly for beams containing soft inclusions, computations were performed using interpolation orders~\(m=2\), \(4\), \(6\), and~\(8\). 
\Cref{tab:OrderInp} summarises the critical strain and critical stress predicted for the unimodal buckling instability for both homogeneous beams (\(E_s=1\)) and soft-inclusion beam~(\(E_s=0.01\)) at various \patRats~\cref{Epatrat} with various orders~\(m\).
For all interpolation orders, the full-domain computations (\(r=1\)) successfully capture the buckling instability and provide consistent predictions of the critical strain and stress---independent of~\(m\) since the interpolation reduces to a direct copy in full-domain computations. 
The patch predictions using \patRat\ \(r<1\) exhibits a dependence on the interpolation order. 
For~\(m=2\) and~\(4\), no buckling instability is detected for either~\(r=0.25\) or~\(r=0.5\), irrespective of the inclusion stiffness. 
These results indicate that low-order interpolation is insufficient to accurately transfer the macroscale information required to reproduce the instability.
In contrast, interpolation orders~\(m=6\) and~\(8\) consistently predict the onset of buckling and yield critical strains and stresses that are in good agreement with the corresponding full-domain results. 
Among these, \(m=8\) provides the most accurate predictions, exhibiting excellent agreement with the full-domain values for both \patRats, with the improvement being particularly pronounced for~\(r=0.25\).
All of the~\(r<1\) computations in~\cref{Sec:NumRes} use interpolation order \(m=8\) to obtain good agreement with \text{the full-domain computations.}

\begin{table}
    \centering
    \caption{\label{tab:OrderInp}
    Critical strain~\(\varepsilon_{\Cr}\) and critical stress~\(\sigma_{\Cr}\) for unimodal buckling instability computed using various orders of the patch interpolation \text{and \patRats.}   }
\begin{equation*}
    \begin{array}{|cc|cc|cc|cc|}
     \hline
&&\multicolumn{2}{c|}{r=0.25}      & \multicolumn{2}{c|}{r=0.5}     & \multicolumn{2}{c|}{r=1} \\
       E_s       & m & \varepsilon_\Cr  & \sigma_\Cr \cdot10^{3}& \varepsilon_\Cr  & \sigma_\Cr \cdot10^{3}& \varepsilon_\Cr  & \sigma_\Cr \cdot10^{3} \\
      \hline
       \multirow{3}{*}{1}                 & \text{2}  & \multicolumn{2}{c|}{\text{no buckling}}  & \multicolumn{2}{c|}{\text{no buckling}} & \multirow{4}{*}{0.093\%} & \multirow{4}{*}{1.022}  \\
                                                    & \text{4}  &  \multicolumn{2}{c|}{\text{no buckling}}  &  \multicolumn{2}{c|}{\text{no buckling}}  &       &   \\
                                                    & \text{6}  & 0.099\%  & 1.104 & 0.093\% & 1.023 &  &  \\
                                                    & \text{8}  & 0.093\%  & 1.024 & 0.093\% & 1.023 &  &  \\
      \hline
      \multirow{3}{*}{\(10^{-2}\)}  & \text{2}  & \multicolumn{2}{c|}{\text{no buckling}} & \multicolumn{2}{c|}{\text{no buckling}} & \multirow{4}{*}{0.079\%} & \multirow{4}{*}{0.493} \\
                                                   & \text{4} & \multicolumn{2}{c|}{\text{no buckling}}  & \multicolumn{2}{c|}{\text{no buckling}}  &                                         &                                   \\
                                                   & \text{6}  & 0.081\% & 0.506 & 0.079\% & 0.493 &     &                                    \\   
                                                   & \text{8}  & 0.078\% & 0.468 & 0.079\% & 0.493 &    &   \\    
     \hline      
\end{array}\end{equation*}
\end{table}%

\subsection{Computing time}
\label{SSct}

As well as having a highly controllable accuracy (\cref{subSec:Emodes,SSbb,SSoi}), a second key advantage of the patch scheme is its computational efficiency---an efficiency which arises from solving the microscale model only within small patches of the system domain. 
Here this efficiency is assessed in two system-level computations: the beam equilibrium calculation and the eigen-mode calculation. 
The first of these encompasses the computational time required to solve the microscale model within the patches, including applying patch coupling via the patch interpolation, until the quasi-Newton iteration reaches the equilibrium solution. 
The latter denotes the computational time required to compute the eigenvalues and eigenvectors of the Jacobian matrix \text{of the patch system.}

\begin{figure}
\centering
\caption{\label{fig:CompTime}
Computational times for the soft-inclusion beam with Young's moduli of the soft inclusions, \(E_s=1\), \(0.1\), \(0.01\), and~\(0.001\), computed in the full domain~(\(r=1\)),  half domain~(\(r=0.5\)), and quarter-domain~\(r=0.25\): (left) beam-equilibrium; (right) eigen-mode.}
\includegraphics[]{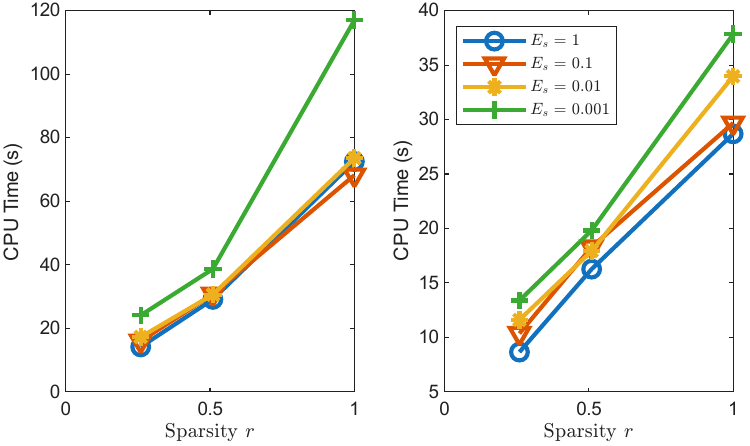}   
\end{figure}

All computations in~\cref{Sec:NumRes} were performed on a standard laptop.\footnote{64-bit Windows, Intel(R) Core(TM) i7-8550U CPU at 1.80\,GHz, 16.0\,GB RAM}
The reported computational times should therefore be interpreted as indicative of the \emph{trends} for the types of system-tasks  in this exploration. 
\Cref{fig:CompTime} shows the computational time required for each case of applied compressive strain.
Typically, an applied strain increment of~\(\delta\varepsilon^e=0.02\%\) then required about forty times as much time to cover \text{the compressive domain~\(0\leq\varepsilon^e\leq0.4\%\). }

In the eigen-mode computation, each applied strain required the computation of the eigenvalues of smallest magnitude---typically the 60~eigenvalues with magnitude less than about~\(0.2\)---together with their corresponding eigenvectors. 
\Cref{fig:CompTime} shows that the computational time decreases approximately linearly with the \patRat~\(r\). 
As Young's modulus~\(E_s\) of the inclusions decreases from~\(1\) to~\(10^{-3}\), the computational time increases slightly, due to the greater complexity of the problem and a consequent reduction in the convergence rate. 
Nevertheless, the reduction in computational time with decreasing \patRat\ remains consistently semi-linear across all investigated values of~\(E_s\). 
These results demonstrate that, within the patch scheme, reducing the computational domain \patRat~\(r\) leads to a nearly proportional reduction in computational time. 
However, the number of patches (which here is proportional to~\(r\)) must remain sufficiently large, in order for the patch spacing~\(H\) to resolve the emergent macroscale \text{structures of the beam.}

\section{Conclusions}
\label{Sec:Cons}

When subjected to a critical compressive load, slender structures may undergo a sudden and significant change in deformation, accompanied by a marked loss of load-carrying capacity and, in severe cases, catastrophic structural failure. 
This phenomenon is referred to as buckling. 
Buckling is inherently multiscale and highly nonlinear, and requires not only resolving the strongly nonlinear structural response in the vicinity of the critical point but also accounting for microscale features that may initiate or alter the instability.
Building on our previous multiscale patch framework for linear elasticity~\citep{TranD2024, TranD2025}, the present work incorporates nonlinear elasticity and the efficient patch scheme to capture large deformations, buckling, and post-buckling behaviour of soft-inclusion beams.
The patch scheme~(\cref{Sec:ProbDef}) is designed to resolve the underlying material microstructure in sufficient detail while maintaining computational efficiency and accuracy in predicting the macroscale response of the beams.
\cite{Maclean2021} provides a freely available Toolbox to empower everybody to try \text{the patch scheme.}

Although this study considers beams containing periodically distributed soft inclusions, the patch scheme is not restricted to this particular microstructure. 
Instead, it provides a general multiscale computational methodology that can be readily applied to a broad class of heterogeneous materials by redefining the microscale discretisation within the patches to represent the underlying material architecture. 
The patch scheme has also been validated to apply to functionally graded micro-heterogeneous material~\cite[e.g.,][]{TranD2025}.
Consequently, the patch scheme serves as a flexible computational framework that can be adapted to different material systems and physical problems \text{with minimal modification.}

%\todo{Are the following predictions cognate to others in the literature?  If so, cite one or two examples here.}
The computational results show that decreasing the stiffness of the soft inclusions lowers both the critical buckling strain and the critical buckling stress, indicating an increased susceptibility of the beam to buckling.
This trend is consistent with previous studies showing that soft inclusions weaken the global mechanical response of heterogeneous beams~\citep{Muller2022} and that reduced inclusion stiffness can lower the buckling resistance of heterogeneous structures~\citep{ZengOlsson2002, HasratiJain2026}.
Soft inclusions have also been shown to substantially modify nonlinear compressive behaviour through localised deformation and buckling mechanisms \citep{SztefekOlsson2009}.
In the present beam, decreasing~\(E_s\) produces the distinctive combination of an earlier buckling onset and a milder global post-buckling response, characterised by smaller stress drops and reduced transverse deflections (\cref{SSbb}).
This behaviour is associated with increasing localisation of deformation within the soft inclusions, which accommodate part of the imposed compression locally and thereby mitigate the \text{global post-buckling deformation.}

The eigenvalues of the Jacobian matrix confirms that buckling is governed by long-wavelength bending modes--the compressive modes remain stable throughout the domain of applied strain explored herein. 
This behaviour is consistent with classical beam buckling theory, in
which the primary instability of a slender beam under axial compression
is associated with its lowest-order bending mode~\citep{TimoshenkoGere1961}.
Patch computations predict the first three critical points of the buckling modes to good accuracy, particularly for the primary unimode, even when only a quarter (\(r=0.25\)) of the beam domain is simulated.
Patch computations also capture the post-buckling behaviour of soft-inclusion beams, although for~\(r=0.25\) there is overestimation of the stress and an underestimation of \text{the beam deflection.}

The accuracy of the patch scheme is controllably influenced by the patch spacing~\(H\) and the interpolation order~\(m\) used to couple neighbouring patches (\cref{subSec:Emodes,SSoi}). 
Low-order interpolations, such as~\(m=2\) or~\(4\), fail to reproduce the buckling instability in reduced-domain computations, indicating insufficient transfer of macroscale information across unsimulated space between patches. 
Higher-order inter-patch interpolation, \(m\ge 6\), accurately captures both the critical strain and critical stress of the buckling instability (users need only change one parameter to the Toolbox by \cite{Maclean2021}).
In particular, order~\(m=8\) provides excellent agreement \text{with full-domain results.}

The patch scheme exhibits an approximately linear reduction in computing time with \patRat~\(r\) (\cref{SSct}). 
Simulations performed with~\(r=0.5\) and~\(r=0.25\) preserve the essential macroscale dynamics while significantly reducing computational time in both equilibrium and eigenvalue calculations.
These results confirm that sparse microscale simulations, when coupled through sufficiently accurate interpolation, provide an efficient and reliable framework for predicting instability phenomena in heterogeneous \text{and microstructured structures.}

\end{document}

%% file: preamble.tex
\IfFileExists{ajrxxx.sty}{% AJRs preference for drafting
\usepackage[a5paper,margin=10mm,tmargin=15mm]{geometry} 
}{% common alternative
}

\usepackage{amsmath,amssymb,microtype,enumitem}
\usepackage[defaultlines=3,all]{nowidow}
\usepackage{pdflscape}
\usepackage{multirow}
\usepackage[leftcaption]{sidecap}
\usepackage{mathtools}
\usepackage{xcolor}
\usepackage[backend=bibtex8 % omit if biber works for you
    ,style=authoryear % for Harvard-like
    ,backref=true % optional
    ]{biblatex}
\IfFileExists{ajr.sty}
  {\bibliography{ajr,bib,bibexport}}
  {\bibliography{bibexport}}
\AtEndDocument{{\raggedright\printbibliography}}
\let\citet\textcite
\let\citep\parencite
\makeatletter
\def\cite{\@ifnextchar[{\parencite}{\textcite}}% only Harvard-like
\makeatother
\DeclareFieldFormat{url}{%
  \iffieldundef{doi}{%
  \url{#1}}{}}
\DeclareFieldFormat{urldate}{%
  \iffieldundef{doi}{%
  \mkbibparens{\bibstring{urlseen}\space#1}}{}}

\usepackage{tikz,pgfplots}
\usetikzlibrary{shapes.geometric}
\pgfplotsset{compat=newest}
\def\tikzsetnextfilename#1{} % define if not using following
\usepgfplotslibrary{external}
\usepackage{doi}
\hypersetup{
    colorlinks   = true, %Colours links instead of ugly boxes
    urlcolor     = blue, %Colour for external hyperlinks
    linkcolor    = blue, %Colour of internal links
    citecolor   = magenta,
    } 
\usepackage[nameinlink,capitalise,noabbrev]{cleveref}
\crefname{equation}{}{}
\let\ref\cref
\let\eqref\cref
\let\autoref\cref

\usepackage[colorinlistoftodos
    ,textsize=footnotesize,textwidth=9em]{todonotes}
\IfFileExists{ajrxxx.sty}{\setuptodonotes{inline,color=orange!20}}
    {\setuptodonotes{color=orange!20}}% common alternative 
  \ifx\undefined
     
  \else
     
  \fi

\def\RaisedName#1{}% to cancel it
\def\ajrSplit#1#2\ajrEndSplit{\def\ajrcar{#1}\def\ajrcdr{#2}}
\def\Vec{}
\renewcommand{\Vec}[1]{%
  \typeout{**** Command #1v denotes vector #1}%
  \expandafter\ajrSplit#1\ajrEndSplit%
  \ifx\ajrcdr\empty % single character
    \expandafter\DeclareRobustCommand\csname#1v\endcsname%
    {{\RaisedName{\ensuremath{\backslash}#1v}% put name above
      \ensuremath{\ifx#1i \vec\imath
        \else\ifx#1j \vec\jmath
        \else\vec #1\fi\fi}% end of typeset it
      \global\csname UsedVec#1true\endcsname% record when used
    }}%
  \else % multicharacter e.g. greek
    \expandafter\DeclareRobustCommand\csname#1v\endcsname%
    {{\ensuremath{\vec{\csname#1\endcsname}}}% record when used
      \global\csname UsedVec#1true\endcsname}%
  \fi%
  \expandafter\newif\csname ifUsedVec#1\endcsname
  \csname UsedVec#1false\endcsname
  \AtEndDocument{\csname ifUsedVec#1\endcsname\else
  \typeout{**** Vec symbol #1v not used.}
  \fi}
  }
\newcommand{\Cal}[1]{% define symbol c#1
  \typeout{**** Command c#1 denotes calligraphic #1}%
  \expandafter\DeclareRobustCommand\csname c#1\endcsname%
  {{\RaisedName{\detokenize{\c}\hspace{-0.5em}#1}% put name above
  \ensuremath{\mathcal #1}% typeset it
  \global\csname UsedCal#1true\endcsname  % record when used
  }}%end defn of \c#1
  \expandafter\newif\csname ifUsedCal#1\endcsname
  \csname UsedCal#1false\endcsname
  \AtEndDocument{\csname ifUsedCal#1\endcsname\else
  \typeout{**** Cal symbol c#1 not used.}
  \fi}
  }

\renewcommand{\vec}[1]{\text{\boldmath$#1$}}

\newcommand{\E}{\cdot 10^}
\def\phalf{^+} % for microgrid index shifting
\def\mhalf{^-} % for microgrid index shifting

\def\pde{\textsc{pde}}

\Vec e\Vec f\Vec k\Vec n\Vec q\Vec x\Vec u
\Vec{varepsilon}\Vec{sigma}\Vec I\Vec J\Vec F\Vec U
\Cal U

\newcommand{\D}[2]{\mathchoice
  {\frac{\partial #2}{\partial #1}}% display
  {{\partial #2}/{\partial #1}}% text
  {{\partial #2}/{\partial #1}}% script
  {{\partial #2}/{\partial #1}}% scriptscript
  }

\newcommand{\Ord}[1]{\ensuremath{\mathcal O%
  \mathchoice{\big(#1\big)}{\big(#1\big)}{(#1)}{(#1)}%
  }}

\newcommand{\opn}[1]{\operatorname{#1}}
\def\E#1{\textsc{e}\ifx#1-{-}\else\ifx#1+{+}\else#1\fi\fi}

\def\uc{blue}
\def\vc{red}
\def\wc{green!60!black}
\def\temp#1{#1}% by default do nothing
\def\oSym{\temp{$\color{\wc}\circledcirc$}}
\def\xSym{\temp{$\color{\wc}\otimes$}}
\def\uSym{\temp{$\color{\uc}\blacktriangleright$}}
\def\vSym{\temp{$\color{\vc}\blacktriangle$}}

\newcommand{\Cr}{\text{cr}}

%% file: Figs/BeamCoord.tex
% Requires: \usepackage{tikz}
% Optional (for nicer arrow heads): \usetikzlibrary{arrows.meta}

\begin{tikzpicture}[scale=1.0]
  % --- Parameters ---
  \def\L{8}     % beam length
  \def\W{2}     % beam height
  \def\r{0.4}   % radius of inclusions

  % --- Beam ---
  \draw[thick] (0,{-\W/2}) rectangle (\L,{\W/2});

  % --- Coordinate axes at origin ---
  \draw[-,thick] (0,0) -- (\L/2-0.5,0);
  \draw[dashed,thick] (\L/2-0.4,0) -- (\L/2+0.3,0);
  \draw[->,thick] (\L/2+0.3,0)  -- (\L+0.5,0) node[below right] {$x$};
  \draw[->,thick] (0,0) -- (0,1.5) node[above left] {$y$};

  % --- Origin mark/label ---
  \fill (0,0) circle (1.2pt);
  \node[below left] at (0,0) {$O$};

  % --- Labels ---
  \node[above left] at (\L+0.5,0) {$L$};
  \node[left] at (-0.1,{\W/2}) {$W/2$};
  \node[left] at (-0.1,{-\W/2}) {$-W/2$};

  % --- Layer of soft inclusions ---
  \foreach \x in {0.75,2.5,5.5,7.25} {
    \draw[fill=gray!30] (\x,0) circle (\r);
  }

\end{tikzpicture}

%% file: Figs/PatchLayout.tex
\def\N{7}\def\Nm{9}\def\Nmm{8} % in x-dirn, N=nx-1, Nm=nx
\def\dM{3}
\def\r{0.2}
\begin{tikzpicture}[x=8mm,y=8mm]
% uniform
% the beam domain

\draw ({-\N+0.5},0) rectangle ({\N+0.5},1);
\foreach \i in {-\N,...,-3} { 
     \draw[dashed] ({\i+1.5},0) -- ({\i+1.5},1);
     \fill[gray] (\i+1,0.5) circle ({\r});
     \draw[dashed] ({-\i-0.5},0) -- ({-\i-0.5},1);
     \fill[gray] (-\i,0.5) circle ({\r});
}
\node at (0.5,0.5) {. . .};
\node at (-\N+1,1.4) {1};
\node at (-\N+2,1.4) {2};
\node at (-\N+3,1.4) {3};
\node at (\N-2,1.4) {$I_{L-2}$};
\node at (\N-1,1.4) {$I_{L-1}$};
\node at (\N,1.4) {$I_L$};

\node[anchor=west] at (-\N,2) {(a) \(r= 1\) , \(I_p=I_L\)};

% 51
\draw[dashed,black] (-\N+0.5,-\dM) rectangle (\N+0.5,-\dM+1);
\draw (-\N+0.5,-\dM) rectangle (-\N+1.5,-\dM+1);
\fill[gray] (-\N+1,-\dM+0.5) circle ({\r});
\draw (-\N+2.5,-\dM) rectangle (-\N+3.5,-\dM+1);
\fill[gray] (-\N+3,-\dM+0.5) circle ({\r});
\draw (-\N+4.5,-\dM) rectangle (-\N+5.5,-\dM+1);
\fill[gray] (-\N+5,-\dM+0.5) circle ({\r});

\draw (\N-4.5,-\dM) rectangle (\N-3.5,-\dM+1);
\fill[gray] (\N-4,-\dM+0.5) circle ({\r});
\draw (\N-2.5,-\dM) rectangle (\N-1.5,-\dM+1);
\fill[gray] (\N-2,-\dM+0.5) circle ({\r});
\draw (\N-0.5,-\dM) rectangle (\N+0.5,-\dM+1);
\fill[gray] (\N,-\dM+0.5) circle ({\r});
\node at (0.5,-\dM+0.5) {. . .};
\node at (-\N+1,-\dM+1.4) {1};
\node at (-\N+3,-\dM+1.4) {2};
\node at (-\N+5,-\dM+1.4) {3};
\node at (\N-4,-\dM+1.4) {$I_{p-2}$};
\node at (\N-2,-\dM+1.4) {$I_{p-1}$};
\node at (\N,-\dM+1.4) {$I_p$};
\node[anchor=west] at (-\N,-\dM+2) {(b) \(r= 0.5\), $I_p=\left(I_{L}+1\right)/2$};

% 26
\draw[dashed,black] (-\N+0.5,-2*\dM) rectangle (\N+0.5,-2*\dM+1);
\draw (-\N+0.5,-2*\dM) rectangle (-\N+1.5,-2*\dM+1);
\fill[gray] (-\N+1,-2*\dM+0.5) circle ({\r});
\draw (-\N+4.5,-2*\dM) rectangle (-\N+5.5,-2*\dM+1);
\fill[gray] (-\N+5,-2*\dM+0.5) circle ({\r});
\draw (\N-4.5,-2*\dM) rectangle (\N-3.5,-2*\dM+1);
\fill[gray] (\N-4,-2*\dM+0.5) circle ({\r});
\draw (\N-0.5,-2*\dM) rectangle (\N+0.5,-2*\dM+1);
\fill[gray] (\N,-2*\dM+0.5) circle ({\r});
\node at (0.5,-2*\dM+0.5) {. . .};
\node at (-\N+1,-2*\dM+1.4) {1};
\node at (-\N+5,-2*\dM+1.4) {2};
\node at (\N-4,-2*\dM+1.4) {$I_{p-1}$};
\node at (\N,-2*\dM+1.4) {$I_p$};
\node[anchor=west] at (-\N,-2*\dM+2) {(c) \(r= 0.25\), $I_p=\left(I_{L}+3\right)/4$};

\end{tikzpicture}

%% file: Figs/Grid.tex
\tikzsetnextfilename{Figs/figpatchgridv}
%\def1{1.6} % scaling of lattice
\def\N{7}\def\Nm{6} % in x-dirn, N=nx-1, Nm=nx
\def\M{6}\def\Mm{5} % in y-dirn, M=ny, Mm=ny-1
\def\r{1.4}
\begin{tikzpicture}[x=14mm,y=14mm]% use these to scale plot coordinates
\filldraw[fill=yellow!15,draw=yellow!15](1.5,1)--(1.5,\M)--(\Nm-0.5,\M)--(\Nm-0.5,1);
\fill[lightgray] (\N/2-0.5,\M/2+0.5) circle (\r);
\foreach \j in {1,...,\Mm} {
       \filldraw[fill=yellow!15,draw=yellow!15](1.5,\j)--(0.5,\j)--(1,\j+0.5)--(0.5,\j+1)--(1.5,\j+1);
        \filldraw[fill=yellow!15,draw=yellow!15](\Nm-0.5,\j)--(\Nm,\j+0.5)--(\Nm-0.5,\j+1);
        \draw[yellow,thick](0.5,\j)--(1,\j+0.5)--(0.5,\j+1);
        \draw[yellow,thick](\Nm-0.5,\j)--(\Nm,\j+0.5)--(\Nm-0.5,\j+1);
      }
\draw [step=1,magenta!40,thin] (0,1) grid (\N-0.3,\M);
\foreach \i in {1,...,\N}   \foreach \j in {1,...,\Mm} % u
      {
        \node[] at (\i-0.5,\j+0.5) {\uSym};           
      }
\foreach \i in {1,...,\N}   \foreach \j in {1,...,\M} % v
      {
        \node[] at (\i-1,\j) {\vSym};     
      }  
\foreach \i in {1,...,\Nm}   \foreach \j in {1,...,\Mm} % sig_xx sig_yy
      {
        \node[] at (\i,\j+0.5) {\oSym};           
      }
\foreach \i in {1,...,\Nm}   \foreach \j in {1,...,\M} % sig_xy
      {
        \node[] at (\i-0.5,\j) {\xSym};     
      }                
\foreach \j in {1,...,\Mm} % u edge fields
    {
     \node[] at (0.5,\j+0.5) {$\square$};    
     \node[] at (\N-0.5,\j+0.5) {$\square$};    
    }          
\foreach \j in {1,...,\M} % v edge fields
    {
     \node[] at (0,\j) {$\square$};    
     \node[] at (\Nm,\j) {$\square$};    
    }      

\node[below] at (-0.1,0.7) {$i=1$}; 
\node[below] at (1,0.7) {$2$}; 
\node[below] at (\N-2,0.7) {$n_x-1$};     
\node[below] at (\N-1,0.7) {$n_x$};     
\node[above] at (-0.5,1-0.2) {$j=1$};      
\node[above] at (-0.5,2-0.2) {$2$};      
\node[above] at (-0.5,\M-1.3) {$n_y-1$}; 
\node[above] at (-0.5,\M-0.2) {$n_y$}; 
\end{tikzpicture}